\documentclass[ejs,preprint,twoside,linksfromyear,nameyear,natbib]{imsart}
\usepackage[colorlinks,citecolor=blue,urlcolor=blue]{hyperref}

\usepackage{mathbh}

\usepackage{amssymb}
\usepackage{latexsym}
\usepackage[centertags]{amsmath}
\usepackage{amsthm}
\usepackage{graphicx, xcolor}
\usepackage{epsfig, subfigure}

\usepackage{array}
\usepackage{algorithmic, algorithm}
\RequirePackage{hypernat}

\doi{10.1214/13-EJS844}
\pubyear{2013}
\volume{7}
\issue{0}
\firstpage{2344}
\lastpage{2371}
\startlocaldefs
\newcommand{\E}{\ensuremath{\mathrm{E}}}
\newcommand{\ud}{\ensuremath{\mathbf{d}}}
\newcommand{\R}{\ensuremath{\mathbf{R}}}
\newcommand{\1}{\ensuremath{\mathbh{1}}}

\DeclareMathOperator{\pen}{\mathbf{pen}}
\DeclareMathOperator{\crit}{\mathbf{crit}}
\DeclareMathOperator{\h}{\mathbf{h}}
\DeclareMathOperator{\KL}{\mathbf{KL}}
\DeclareMathOperator{\BIC}{\mathbf{BIC}}
\DeclareMathOperator{\AIC}{\mathbf{AIC}}
\DeclareMathOperator{\ICL}{\mathbf{ICL}}

\DeclareMathOperator{\CteDim}{\mathbf{Cte*Dim}}
\newcommand{\argmin}{\mathop{\mathrm{arg\,min}}}
\newcommand{\argmax}{\mathop{\mathrm{arg\,max}}}

\newcommand{\Ss}{\ensuremath{\mathbb{S}}}
\newcommand{\N}{\ensuremath{\mathbb{N}}}
\newcommand{\LL}{\ensuremath{\mathbb{L}}}

\newcommand{\colmodels}{\ensuremath{\mathcal{C}}} %\mathbb{M}

\newtheorem{thm}{Theorem}
\newtheorem{lem}{Lemma}

\newtheorem{prop}{Proposition}

\theoremstyle{definition}

\allowdisplaybreaks

\endlocaldefs

\begin{document}
\begin{frontmatter}

\title{Clustering and variable selection for categorical multivariate data}
\runtitle{Clustering and variable selection for categorical data}

\author{\fnms{Dominique} \snm{Bontemps}\thanksref{t1}\corref{}\ead[label=e1]{dominique.bontemps@math.univ-toulouse.fr}}
\address{\printead{e1}}
\thankstext{t1}{Supported by Institut de Math\'ematiques de Toulouse, Universit\'e de Toulouse}
\medskip
\and
\author{\fnms{Wilson} \snm{Toussile}\thanksref{t2}\ead[label=e2]{wilson.toussile@u-psud.fr}}
\address{\printead{e2}}
\thankstext{t2}{Supported by INSERM U669 and Fac. de M\'edecine, Universit\'e Paris-Sud 11}

\runauthor{D. Bontemps and W. Toussile}

\begin{abstract}
This article investigates unsupervised classification techniques for
categorical multivariate data. The study employs multivariate
multinomial mixture modeling, which is a type of model particularly
applicable to multilocus genotypic data. A model selection procedure is
used to simultaneously select the number of components and the relevant
variables. A~non-asymptotic oracle inequality is obtained, leading to
the proposal of a new penalized maximum likelihood criterion. The
selected model proves to be asymptotically consistent under weak
assumptions on the true probability underlying the observations. The
main theoretical result obtained in this study suggests a penalty
function defined to within a multiplicative parameter. In practice, the
data-driven calibration of the penalty function is made possible by
slope heuristics. Based on simulated data, this procedure is found to
improve the performance of the selection procedure with respect to
classical criteria such as $\mathbf{BIC}$ and $\mathbf{AIC}$. The new
criterion provides an answer to the question ``Which criterion for
which sample size?'' Examples of real dataset applications are also
provided.
\end{abstract}


% history:
\received{\smonth{3} \syear{2012}}

\end{frontmatter}

%s1 ###
\section{Introduction} \label{Sect:Introduction}

This article investigates unsupervised classification and variable
selection in the context of categorical multivariate data. Considering
the frequencies of each variable's categories, the underlying
population is assumed to be structured into sub-populations of a
certain unknown number $K$. The possibility exists that only a subset
$S$ of the variables are relevant for clustering purposes. This subset
$S$ may significantly influence the interpretation of results.

Building on \cite{ToussileGassiat2009}, we consider the modeling
problem of simultaneously selecting $K$ and $S$ in a density
estimation framework. A penalized maximum likelihood procedure is used,
which also permits us to estimate the frequencies of the categories at
the same time. Individuals are subsequently clustered with the maximum
a posteriori (MAP) method. Our study offers a data-driven model
selection criterion derived from a new non-asymptotic oracle
inequality.

Clustering of categorical multivariate data is %an important topic
used in many fields such as social sciences, health, marketing,
population genetics, etc. \cite[see for instance][]{McCutcheon,
McLachlanPeel, CollinsLanza}. The population genetics specific
framework we develop applies to multilocus genotypic data. This type of
data corresponds to a situation in which each variable describes the
generic variants or \emph{alleles} of a genetic marker (called a
\emph{locus}). For diploid organisms, which have one allele from each
of their parents, two unordered alleles are observed at each locus.

We use finite mixture models to investigate clustering in discrete
settings, under the common hypothesis that the variables are
conditionally independent with respect to each component of the
mixture. In the literature, such models are also known as latent class
models, which were first introduced by \cite{Goodman74}. The family of
latent class models has proven to be successful in many practical
situations \citep[see for instance][]{Rigouste_Cappe_Yvon_2006}.

\begin{sloppypar}
Various model-based clustering methods for categorical multivariate
data have been proposed in recent years
\citep{CeleuxGovaert91,Pritchard2000, Chen2006, Corander2008}. Several
of these papers used a Bayesian approach \cite[for details,
see][]{Celeux99computational, Rigouste_Cappe_Yvon_2006}. Yet, the
problem of variable selection in clustering for categorical
multivariate data was first addressed in \cite{ToussileGassiat2009}.
The simulated data used in their study suggested that a
variable selection procedure %based on the Bayesian Information
%Criterion ($\BIC$)
could significantly improve clustering and prediction capacities for
our intended
framework. Furthermore, the article provided theoretical consistency
results for $\BIC$ type criteria. Such criteria are, however, known
to require large sample sizes to attain their consistency behavior
in discrete settings \cite[see also][]{Nadif1998}.
\end{sloppypar}

We adopt an oracle approach to conduct the present study. It is not our
aim to choose the true model $\mathcal{M}_{(K_0,\ S_0)}$ underlying the
data, although our procedure is found to also perform well in that
respect. Instead, it is our intention to propose a criterion that is
designed to minimize a risk function based on the Kullback-Leibler
divergence of the estimated density with respect to the true density.
%some risk function of the estimated density with respect to the true
%density.
In this context, ``simpler'' models are preferable to
$\mathcal{M}_{(K_0,\ S_0)}$, for which too many parameters may result
in estimators that over fit the data. In fact, it is unnecessary to
assume that $P_0$ belongs to one of the competing
models~$\mathcal{M}_{(K,S)}$.

The non-asymptotic penalized criterion we propose in this paper is
based on the metric entropy theory and a theorem of \cite{Massart2007}.
The new criterion leads to a non-asymptotic oracle inequality, which
compares the risk of the selected estimator with the risk of the estimator
that is associated with the (unknown) best model (see Theorem~\ref
{thm:main-theorem}
below). A large volume of literature examines model selection through
penalization from a non-asymptotic perspective. Research in this area
is still in development and follows the emergence of new sophisticated
tools of probability, such as concentration and deviation inequalities
\citep[see][and references therein]{Massart2007}. This kind of
approach has only recently been applied to mixture models; \cite
{MaugisMichel2008thm}
were the first to use it for Gaussian mixture models. Our study focused
on discrete variables.

\begin{sloppypar}
Nevertheless, the obtained penalty function presents certain drawbacks:
The function depends on a multiplicative constant for which sharp upper
bounds are not available, and it leads, in practice, to an over
penalization that is even worse than $\BIC$. We therefore calibrate the
constant with the so-called slope heuristics proposed in
\cite{BirgeMassart2007}. Slope heuristics, although only fully
theoretically validated in the Gaussian homoscedastic and
heteroscedastic regression frameworks \citep{BirgeMassart2007,
ArlotMassart2009}, have been implemented in several other frameworks
\citep[see][for applications in density estimation, genomics,
etc.]{MaugisMichel2008slope,LebarbierPhDThesis,Verzelen2009,Villers2007}.
The simulations described in
Subsection~\ref{sect:Simulated_experiments} illustrate that our
criterion behaves well with respect to more classical criteria such as
$\BIC$ and $\AIC$, both in terms of density estimation (even when $n$
is relatively small) and true model selection. The criterion can be
considered part of the family of General Information Criteria
\citep[see for instance][whose criterion presents some analogy to slope
heuristics]{BaiRaoWu1999}.
\end{sloppypar}

Section~\ref{Sect:Model_methods} of this paper presents the mixture
model framework and the model selection paradigm. In
Section~\ref{Sect:Main_result} we describe and prove our main result,
the oracle inequality. Section~\ref{Sect:In_practice} examines the
practical aspect of our procedure, which was implemented in the
stand-alone software \texttt{MixMoGenD} (Mixture Model using Genotypic
Data) that was first introduced in \cite{ToussileGassiat2009}.
Simulated experimental results are presented in
Section~\ref{sect:Simulated_experiments}, including a comparison of our
proposed criterion with classical $\BIC$ and $\AIC$, considering both
the selection of the true model and the density estimation. Examples
of applications to real datasets can be found in
Section~\ref{sect:real}. Finally, the Appendices contain several
technical results used in the main analysis.

%s2 ###
\section{Models and methods}\label{Sect:Model_methods}

%s2.1 ###
\subsection{Framework}\label{Subsect:Framework}

Consider independent and identically distributed (iid) instances
of a multivariate random vector $X=(X^{l})_{1\leq l\leq L}$,
where the number of categorical variables $L$ is potentially large.
We investigate two main settings:
\begin{enumerate}
\item\label{config:simple-multinomial} Each $X^l$ is a multinomial
    variable taking values in $\{1,\ldots,\ A_l\}$.
\item\label{config:double-multinomial} Each $X^l$ consists of a
    (unordered) set $\{X^{l,1},\ X^{l,2}\}$ of two (possibly equal)
    qualitative variables taking their values in the same set $\{
    1,\ldots,\ A_l\}$.
\end{enumerate}
Throughout this article, these two settings are referred to as
Case~\ref
{config:simple-multinomial}
and Case~\ref{config:double-multinomial}. In both cases, numbers
denoted by $A_l$ are assumed to be known and to satisfy $A_l\geq2$.

Case~\ref{config:simple-multinomial} is generic, whereas Case~\ref
{config:double-multinomial}
is more specific to multilocus data. Our results (presented
below) could easily be extended to other kinds of discrete models,
provided it is possible to compute the metric entropies as described
in Section~\ref{Sect:Main_result}.

The studied sample is assumed to originate from a population structured
into a certain (unknown) number $K$ of sub-populations (clusters),
where each cluster is characterized by a set of category frequencies.
The (unobserved) sub-population an individual comes from is denoted by
the variable $Z$, which takes its values in the set $\{1,\ldots, k,
\ldots,\ K\}$ of the different cluster labels. The distribution of $Z$
is given by the vector $\mathbf{\pi} = (\pi_k)_{1\leq k \leq K}$, where
$\pi_k=P(Z=k)$. Variables $X^1,\ldots,\ X^L$ are assumed to be
conditionally independent given $Z$. For
Case~\ref{config:double-multinomial}, the $X^{l,1}$ and $X^{l,2}$
states of the $l^{\text{th}}$ variable are also assumed to be
conditionally independent given $Z$. In accordance with these
assumptions, the probability distribution of an observation
$x=(x^l)_{1\leq l\leq L}$ in a population $k$ is given in the following
equations:
\begin{align}
P\left(x|\ Z=k\right) &= \prod_{l=1}^L P\left(x^l | Z=k\right)
\notag
\\
\text{Case~\ref{config:simple-multinomial}: } P\left(x^l | Z=k\right)
&= \alpha_{k,l,x^{l}} \notag\\
\text{Case~\ref{config:double-multinomial}: } P\left(x^l | Z=k\right)
&= \left(2-\1_{x^{l,1}=x^{l,2}}\right) \alpha_{k,l,x^{l,1}} \alpha
_{k,l,x^{l,2}} \label{HW_modele}
\end{align}
where $\alpha_{k,l,j}$ is the probability of the modality $j$
associated with the variable $X^l$ in population $k$. The mixing
proportions $\pi_k$ and the probabilities $\alpha_{k,l,j}$ are treated
as parameters.

These assumptions, which are considered classical in latent class
model literature, are known as \emph{Linkage Equilibrium} (LE) and
\emph{Hardy-Weinberg Equilibrium} (HWE) in the context of genomics.
Such assumptions may seem simplistic because they disregard the
migrations between populations and assume that the parents of a given
individual are taken uniformly at random in the population to which the
individual belongs. Nevertheless, these assumptions have proven
useful in describing many population genetic attributes, and they
continue to serve as a base model in the development of more realistic
models of microevolution.

Simplified and misspecified models are often preferable to achieve greater
precision with the oracle approach (as explained in the introduction).
Use of these preferred models introduces a modeling bias in order to
obtain more robust estimators and classifiers and also leads to a
smaller estimation error. In particular, the
introduction of covariances is unlikely to produce better estimates
because it would increase the dimensions of the considered models.
This fact also justifies the following simplification:

It is possible that the structure of interest is contained in only a
subset $S$ of the available variables $L$; the other variables may be
useless and could even hinder the detection of a reasonable clustering
into statistically different populations. The frequencies of the
categories are different in at least two populations for the variables
in $S$; we refer to them as clustering variables. For the other
variables, the categories are assumed to be equally distributed across
the clusters. The simulations performed in \cite{ToussileGassiat2009}
illustrate the benefits of this approximation.

In our case, $\beta_{l,j}$ denotes the frequency of the category $j$
associated with the variable $X^l$ in the whole population:
\begin{equation*}
\beta_{l,j}=\alpha_{1,l,j}=\cdots=\alpha_{k,l,j}=\cdots=\alpha_{K,l,j}\
\textrm{for any } l\notin S \textrm{ and } 1\leq j \leq A_l.
\end{equation*}
Clearly, $S=\emptyset$ if $K=1$, otherwise $S$ belongs to $\mathcal
{P}^{*}(L)$,
which is the set of all nonempty subsets of $\{1,\ldots,L\}$.

Summarizing these assumptions, we can express the likelihood of an
observation $x=(x^l)_{1\leq l\leq L}$:
\begin{align}
\begin{split}
\text{Case~\ref{config:simple-multinomial}: } P_{(K,\, S, \theta)}(x)
&=\left[\sum_{k=1}^K \pi_k \prod_{l\in S} \alpha_{k,l,x^{l}}\right]
\times\prod_{l\notin S} \beta_{l,x^{l}} \\
\text{Case~\ref{config:double-multinomial}: } P_{(K,\, S, \theta)}(x)
&= \left[\sum_{k=1}^K \pi_k \prod_{l\in S}\left(2-\1
_{x^{l,1}=x^{l,2}}\right)
\alpha_{k,l,x^{l,1}}\times\alpha_{k,l,x^{l,2}}\right]\\
&\quad\times\prod_{l\notin S} \left(2-\1_{x^{l,1}=x^{l,2}}\right
)\beta
_{l,x^{l,1}}\beta_{l,x^{l,2}}
\end{split}
\label{Melange_selection}
\end{align}
where $\theta=(\mathbf{\pi},\mathbf{\alpha},\mathbf{\beta })$ is a
multidimensional parameter with
\begin{align*}
\mathbf{\alpha} &=\left(\alpha_{k,l,j}\right)_{1\leq k \leq K;\
l\in
S;\ 1\leq j \leq A_l}\\
\mathbf{\beta} &=\left(\beta_{l,j}\right)_{l\notin S;\ 1\leq j
\leq A_l}.
\end{align*}
For a given $K$ and $S$, $\theta=\theta_{(K,\, S)}$ ranges in the
set
%
%e1 ###
\begin{equation} \label{Theta-KS}
\Theta_{(K,\, S)} = \Ss_{K-1}\times\bigg[\prod_{l\in S} \Ss
_{A_l-1}\bigg]^K \times
\prod_{l\notin S}\mathbb{S}_{A_l-1},
\end{equation}
where $\Ss_{r-1}=\{p=(p_1,\ p_2,\ldots,\ p_r)\in [0,\ 1]^r :
\sum_{j=1}^r p_j = 1\}$ is the $(r-1)$-dimensional simplex.

Then, we consider the collection of all parametric models
%
%e2 ###
\begin{equation} \label{model_KS}
\mathcal{M}_{(K,\, S)}=\left\{P_{(K,\, S, \theta)}: \theta\in
\Theta
_{(K,\, S)} \right\}
\end{equation}
with $(K,\, S) \in\colmodels:= \{(1, \emptyset)\} \cup(\N
\backslash
\{0, 1\}) \times\mathcal{P}^{*}(L)$.
To minimize the use of notations, we make frequent use of the single
index $m\in\colmodels$ instead of using $(K,\, S)$.

Each model $\mathcal{M}_{(K,\, S)}$ corresponds to a particular
structure situation with $K$ clusters and a subset $S$ of clustering
variables. Inferring $K$ and $S$ presents a model selection problem in
a density estimation framework and also leads to data clustering
through the estimation $\widehat{\theta}$ of the parameter
$\theta_{(K,\, S)}$ and the prediction of the class $z$ of an
observation $x$ by the MAP method:
\[
\widehat{z}=\argmax_{1\leq k\leq K}P_{(K,\, S,\widehat{\theta
})}\left
(Z=k|X=x\right).
\]

%s2.2 ###
\subsection{Maximum Likelihood Estimation}

Our study implements the maximum likelihood estimator (MLE). For each
model $\mathcal{M}_{(K,S)}$, the MLE corresponds to the minimum
contrast estimator\break $\widehat{P}_{(K,S)} = P_{(K,\, S,
\widehat{\theta})}$ of the log-likelihood contrast
%
%e3 ###
\begin{equation} \label{eq:log-likelihood-contrast}
\gamma_n\left(P\right)=- \dfrac{1}{n}\sum_{i=1}^n\ln P\left
(X_i\right).
\end{equation}

The Kullback-Leibler divergence $\KL$ provides a suitable risk function
that measures the quality of an estimator in a density estimation.
However, this divergence also has disadvantages in the context of a discrete
framework; in fact, the MLE assigns a zero probability to any unobserved
categories in the sample. Consequently, the Kullback-Leibler risk
%$
%
%e4 ###
\begin{equation} \label{eq:KL-risk}
%R_{\KL}\left(K,\, S\right)=
\E_{P_0}\left[\KL\left(P_0,\ \widehat{P}_{\left(K,\, S\right
)}\right
)\right]
\end{equation}
%
%$
is infinite. In the following, we therefore consider a slightly different
collection $\colmodels^{\varepsilon}$ of competing models $\mathcal
{M}_m^{\varepsilon}$
coupled with a threshold on the parameters of $\varepsilon>0$ :
% in $\alpha$ and $\beta$:
%
\begin{equation*}
\mathcal{M}_{K,S}^{\varepsilon}:=\left\{P_{K,S,\theta}(\cdot)|\
\theta
=(\pi,\alpha,\beta)\in\Theta_{K,S},\ \alpha_{k,l,j}\geq
\varepsilon
\textrm{ and }\beta_{l,j}\geq\varepsilon,\ \forall k,l,j\right\}.
\end{equation*}
Such models disqualify probability distributions that assign too low of
a probability (particularly zero) to certain categories. A good choice
of $\varepsilon$ can result in the same collection of maximum
likelihood estimators with a probability tending to one, as was
demonstrated by a result obtained in \cite[Appendix
D]{ToussileGassiat2009}: if the true probability $P_0$ of the
observations is positive, then for any $(K,\, S)$ a real
$\varepsilon=\varepsilon_{K,S}>0$ exists, such that
\begin{equation*} %\label{theorem2.1}
-\gamma_n\left(\widehat{P}_{\left(K,\, S\right)}\right)=\sup
_{P\in
\mathcal{M}_{K,S}^{\varepsilon}}\left\{-\gamma_n(P)\right\}
+o_{P_0}\left
(1\right),
\end{equation*}
where $\widehat{P}_{(K,\, S)}$ is the MLE in a non-truncated model.

For the sake of simplicity, $\mathcal{M}_{(K,S)}^{\varepsilon}$ is also
denoted by $\mathcal{M}_{(K,S)}$, and $\widehat{P}_{(K,\, S)}$
represents a minimizer of the contrast $\gamma_n$ within $\mathcal
{M}_{(K,S)}^{\varepsilon}$. Additionally, because we cannot discover
more than $n$ clusters from an $n$-sample, we only consider the models
indexed by $(K,S)$ for which the number of clusters $K$ is smaller than
the sample size $n$.

Let $(K^*,\, S^*)$ be a minimizer over $(K,\, S)$ of the
Kullback-Leibler risk (\ref{eq:KL-risk}). The ideal candidate density
$\widehat{P}_{(K^*,\, S^*)}$, or oracle density, is not accessible
because it is dependent on the (unknown) true density $P_0$. The oracle
density is used as a benchmark to quantify the quality of our model
selection procedure: the simulation performed in
paragraph~\ref{sect:oracle_experiment} compares the Kullback-Leibler
risk of the selected estimator $\widehat{P}_{(\widehat{K}_n,\
\widehat{S}_n)}$ with the oracle risk.

%s2.3 ###
\subsection{Model selection through penalization}

\label{Subsect:Model_selection}

Minimization of a penalized contrast is a common method to solve model
selection problems. The selected model $\mathcal{M}_{(\widehat
{K}_{n},\widehat{S}_{n}}$ is a minimizer of a penalized criterion of
the form
\begin{equation*}
\crit(K,S) = \gamma_n\big(\widehat{P}_{(K,S)}\big) + \pen_n(K,S),
\end{equation*}
in which $\pen_n:\ \colmodels\rightarrow\mathbb{R}_+$ is the
penalty function. Eventually, the selected estimator becomes $\widehat
{P}_{(\widehat{K}_n,\ \widehat{S}_n)}$.

The penalty function is designed to avoid over-fit problems. Classical
penalties, such as those used in $\AIC$ and $\BIC$ criteria, are
based on model dimensions. In the following, we refer to the number
of free parameters %of a model $\mathcal{M}_{(K,\, S)}$% given by
%
%e5 ###
\begin{equation} \label{eq:Dimension}
D_{(K,\, S)} = K-1 + K \sum_{l\in S}\left(A_l-1\right) + \sum
_{l\notin
S}\left(A_l-1\right)
\end{equation}
as the dimension of the model $\mathcal{M}_{(K,\, S)}$. The penalty
functions of $\AIC$ and $\BIC$ are respectively defined by
\begin{equation*} %\label{eq:def-BIC-AIC}
\begin{split}
\pen_{\AIC}\left(m\right)&=\frac{1}{n}\, D_m;\\ %\cdot
\pen_{\BIC}\left(m\right)&=\frac{\ln n}{2n}\, D_m.
\end{split}
\end{equation*}

%s3 ###
\section{New criteria and non-asymptotic risk bounds}

\label{Sect:Main_result}

%s3.1 ###
\subsection{Main result}

\label{Subsect:Main_result}

Our main result provides an oracle inequality for Case~\ref
{config:simple-multinomial} and Case~\ref{config:double-multinomial}.
This inequality links the Hellinger risk
$\E_{P_0}[\h^2(P_0,\widehat{P}_{(\widehat {K}_n,\widehat{S}_n)})]$ of
the selected estimator to the Kullback-Leibler divergence $\KL$ between
the true density and each model in the model collection. Recall that
for two probability distributions $P$ and~$Q$, and given $s$ and $t$ as
density functions with respect to a common $\sigma$-finite
measure~$\mu$, the Hellinger distance between $P$ and $Q$ is the
quantity $\h(P,Q)$ defined by %(improperly)
%
%e6 ###
\begin{equation} \label{eq:def-hellinger}
\h(P,Q)^2 = \int\left(\sqrt{s(x)} - \sqrt{t(x)}\right)^2 d\mu
(x). %
%\frac{1}{2}
\end{equation}

Unlike $\KL$, which is not a metric, the Hellinger distance $\h$
allowed us to take advantage of the metric properties (metric entropy)
of the models. %This work is done in Theorem~\ref{thm:Massart}
%presented in the next subsection, which relies on Talagrand's
%inequality. A direct work on this Talagrand's inequality
Use of the metric entropy may be avoided by directly investigating
Talagrand's inequality, which forms the basis of our study. Such
an investigation may then lead to an oracle inequality directly on
the Kullback-Leibler risk (\ref{eq:KL-risk}) and with explicit constants;
this remains to be done in our context. See \cite[S 2.2]{MaugisMichel2008thm}
for more insight regarding this topic.

\begin{thm} \label{thm:main-theorem} We consider the collection
$\colmodels$ of the models defined above and a corresponding collection
of $\rho$-MLEs $\big(\widehat{P}_{(K,S)}\big)_{(K,S)\in\colmodels}$.
Thus, for every \mbox{$(K,\, S)\in\colmodels$}, we obtain
\[
\gamma_n\big(\widehat{P}_{(K,S)}\big) \leq\inf_{Q\in\mathcal
{M}_{(K,\, S)}} \gamma_n(Q) + \rho.
\]

Let $A_{\max}=\sup_{1\leq l\leq L}A_{l}$, and let $\xi$ be defined by
$\xi=\frac{4\sqrt{LA_{\max}}}{2^{L+1}-1}$ in
Case~\ref{config:simple-multinomial} and $\xi=\frac{4\sqrt {LA_{\max
}}}{2^{2L+1}-1}$ in Case~\ref{config:double-multinomial}, and assume
$\xi\leq1$.

There exist absolute constants $\kappa$ and $C$, such that whenever
%
%e7 ###
\begin{equation} \label{eq:pen}
\pen_n(K,S) \geq\kappa\, \Bigg(5+\sqrt{\max\left(\frac{\ln
n+\ln
L}{2},\; \frac{\ln2}{2}+\ln L\right)}\Bigg)^2 \frac{D_{(K,\, S)}}{n}
\end{equation}
for every $(K,\, S)\in{\colmodels}$, then the model $\mathcal {M}_{
(\widehat{K}_n,\widehat{S}_n)}$ exists, where
$(\widehat{K}_n,\widehat{S}_n)$ minimizes
\[
\crit(K,S)= \gamma_n\big(\widehat{P}_{(K,S)}\big)+\pen_n(K,S)
\]
over ${\colmodels}$. Furthermore, whatever the underlying probability
$P_0$,
\begin{multline*}
\E_{P_0}\left[\h^2\left(P_0,\widehat{P}_{\left(\widehat
{K}_n,\widehat
{S}_n\right)}\right)\right] \\
\leq C\, \left( \inf_{(K,\, S)\in\colmodels} \left(\KL\left
(P_0,\mathcal{M}_{(K,\, S)}\right) + \pen_n(K,\, S)\right) + \rho+
\frac{(3/4)^L}{n}\right)
\end{multline*}
where, for every $(K,\, S) \in\colmodels$, $\KL(P_0,\mathcal {M}_{(K,\,
S)}) = \inf_{Q\in\mathcal{M}_{(K,\, S)}} \KL(P_0,Q)$.
\end{thm}
We use the condition $\xi\leq1$ to avoid
more complicated calculations in our proof. In practice, $\xi$ is very likely
to be smaller than $1$ (unless $L$ is very small).

Consider the following:
\begin{itemize}
\item Theorem~\ref{thm:main-theorem} is a density estimation result: it
quantifies the quality of the parameter estimation, which defines
the mixture density components. However, it is not easily connected
to a classification result.
\item The leading term of the penalty for large $n$ is
    $\kappa\tfrac {\ln n}{2}\,\tfrac{D_{(K,\, S)}}{n}$, which is a
    $\BIC$ type penalty function. Consequently, we can apply
    Theorem~2 from \cite{ToussileGassiat2009}: when the underlying
    distribution $P_0$ belongs to one of the competing models, the
    smallest model $(K_0, S_0)$ containing $P_0$ is selected with a
    probability tending to $1$ as $n$ approaches infinity.
\item Such a penalty is not surprising in our context; it is, in fact, very
similar to the penalty obtained by \cite{MaugisMichel2008thm} for a Gaussian
mixture framework.
\item Sharp estimates of $\kappa$ are not available. In practice,
Theorem~\ref{thm:main-theorem}
is too conservative and leads to an over-penalized criterion that
is outperformed by smaller penalties. Therefore, Theorem~\ref
{thm:main-theorem} is mainly used
to suggest the shape of the penalty function
%
%e8 ###
\begin{equation} \label{eq:shape.penalty}
\pen_n(K,\, S) = \lambda\, \frac{D_{(K,\, S)}}{n}
\end{equation}
where the parameter $\lambda$ is chosen depending on $n$ and the
collection $\colmodels$ --- but not on $(K,S)$. Slope heuristics
\citep{BirgeMassart2007,ArlotMassart2009} can be used in practice
to calibrate $\lambda$. This is done in Section~\ref{Sect:In_practice},
where we use change-point detection \cite[see][]{LebarbierPhDThesis}
in connection with slope heuristics.
\item Because $\h^2$ is upper bounded by $2$, the non-asymptotic feature
of Theorem~\ref{thm:main-theorem} becomes important when $n$ is
large enough with respect to $D_{(K,\, S)}$. But even with small
values of $n$, the simulations performed in Subsection~\ref
{sect:Simulated_experiments}
show that the penalized criterion that is calibrated by using slope
heuristics maintains good behavior.
\end{itemize}

%s3.2 ###
\subsection{A general tool for model selection}

\label{Subsect:General_tool}

Theorem~\ref{thm:main-theorem} is obtained from \cite
[Theorem~7.11]{Massart2007},
whose research investigated model selection problems by proposing
penalty functions related to geometrical properties of the models,
namely metric entropy with bracketing for the Hellinger distance.

We examine the following framework: Consider some measurable space
$(A, {\mathcal A})$, and $\mu$ as a $\sigma$-finite positive measure
on $A$. A collection of models $(\mathcal{M}_m)_{m \in\colmodels}$
is given, where each model $\mathcal{M}_{m}$ is a set of probability
density functions $s$ with respect to $\mu$. The following relation
permits us to extend the definition of $\h$ to the positive functions
$s$ or $t$, whose integral is finite but not necessarily $1$. The
function defined by $\sqrt{s}(x) = \sqrt{s(x)}$ is denoted by $\sqrt
{s}$, and $\| \cdot\|_2$ denotes the usual norm in $\LL^2(\mu)$;
then
\begin{equation*}
\h(s,t) = \|\sqrt{s} - \sqrt{t} \|_2.
\end{equation*}

To restate the definition of metric entropy with bracketing, consider
some collection $F$ of measurable functions on $A$ and $d$ as
one of the following metrics on $F$: $\h$, $\|\cdot\|_{1}$, or $\|
\cdot
\|_{2}$.
A bracket $[l,u]$ is the collection of all measurable functions $f$
such that $l\leq f\leq u$. Its $d$-diameter is the distance $d(u,l)$.
Then, for every positive number $\varepsilon$, $N_{[\cdot
]}(\varepsilon,F,d)$
denotes the minimal number of brackets whose $d$-diameter is no
larger than $\varepsilon$, which is required to cover $F$. The
$d$-entropy with bracketing of $F$ is defined as the logarithm of
$N_{[\cdot]}(\varepsilon,F,d)$ and is denoted by $H_{[\cdot
]}(\varepsilon,F,d)$.

We assume that for each model $\mathcal{M}_m$ the square entropy
with bracketing $\sqrt{H_{[\cdot]}(\varepsilon,\mathcal{M}_m,\h)}$
is integrable at $0$. Consider some function $\phi_m$ on $\R_{+}$
with the following properties:

\begin{enumerate}
\setcounter{enumi}{8}
\renewcommand{\theenumi}{\textup{(\Alph{enumi})}}
%\textbf{}
%
\item\label{cond:i} $\phi_m$ is nondecreasing, $x\mapsto\phi_m(x)/x$
is nonincreasing on $(0,+\infty)$ and for every $\sigma\in\R_+$ and
every $u \in\mathcal{M}_m$
\[
\int_{0}^{\sigma} \sqrt{H_{[\cdot]}\left(x, S_m(u,\sigma), \h\right)}
\ud x \leq\phi_m(\sigma),
\]
where $S_m(u,\sigma) = \{ t\in\mathcal{M}_m: \|\sqrt{t} - \sqrt
{u}\|_2 \leq\sigma\}$.
\end{enumerate}

\ref{cond:i} is satisfied, in particular with $\phi_m(\sigma) = \int
_{0}^{\sigma} \sqrt{H_{[\cdot]}(x, \mathcal{M}_m, \h)} \ud
x$.\vadjust{\goodbreak}

%In order to avoid measurability problems, we suppose that for each $m
%\in{\colmodels}$, the following separability condition is satisfied
%for $\mathcal{M}_m$:
%\begin{enumerate}
% \setcounter{enumi}{12}
% \renewcommand{\theenumi}{\textup{(\Alph{enumi})}}%\textbf{}
% \item\label{cond:M} There exists some countable subset $
%\mathcal{M}_m'$ of $\mathcal{M}_m$ and a set $A'\subset A$ with $
%\mu(A')=\mu(A)$ such that for every $t\in\mathcal{M}_m$, there exists
%some sequence $(t_k)_{k\geq1}$ of elements of $\mathcal{M}_m'$ such
%that for every $x\in A'$, $\ln(t_k(x))$ tends to $\ln(t(x))$ as $k$
%tends to infinity.
%\end{enumerate}
\cite{Massart2007} stated a separability
condition, which was denoted (M) in the text, to avoid measurability
problems. This condition is easy
to verify in our context and we omit it for greater legibility of
the theorem.

\begin{thm} \label{thm:Massart} %[\cite[Theorem~7.11]{Massart2007}]
Let $X_1,\ldots,X_n$ be iid random variables with an unknown
density $s$ with respect to some positive measure $\mu$. Let $\{
\mathcal
{M}_m\}_{m \in\colmodels}$
be some at most countable collection of models. %, each fulfilling
%\ref{cond:M}.
We consider a corresponding collection of $\rho$-MLEs
$(\widehat{s}_m)_{m}$. Let $\{x_m\}_{m \in\colmodels}$ be some family
of nonnegative numbers such that\looseness=1
\[
\sum_{m \in\colmodels} e^{-x_m} = \Sigma< \infty,
\]\looseness=0
and for every ${m \in\colmodels}$, considering $\phi_m$ with property
\ref{cond:i}, define $\sigma_m$ as the unique positive solution of the
equation
%
%e9 ###
\begin{equation} \label{eq:def-sigma-m}
\phi_m(\sigma) = \sqrt{n} \sigma^2.
\end{equation}
Let $\pen_n: {\colmodels} \rightarrow\R_+$ and consider the penalized
log-likelihood criterion
\[
\crit(m) = \gamma_n\left(\widehat{s}_m\right) + \pen_n(m).
\]
Then, some absolute constants $\kappa$ and $C$ exist, such that
whenever
\begin{equation*}
\pen_n(m) \geq\kappa\, \left(\sigma_m^2 + \frac{x_m}{n}\right)\
\text
{for every}\ m\in{\colmodels},
\end{equation*}
some random variable $\widehat{m}$ that minimizes $\crit$ over
${\colmodels}$
exists. Furthermore, whatever the density $s$,
\begin{equation*}
\E_s\left[\h^2\left(s,\widehat{s}_{\widehat{m}}\right)\right]
\leq C\,
\left( \inf_{m\in\colmodels} \left(\KL\left(s,\mathcal
{M}_m\right) +
\pen_n(m)\right) + \rho+ \frac{\Sigma}{n}\right).
\end{equation*}
\end{thm}

Concerning Theorem~\ref{thm:Massart}, \cite{Massart2007} explained that
$\sigma_m^2$ has the role of a variance term of $\widehat{s}_m$,
whereas the weights $x_m$ take into account the number of models
$m$ of the same dimension.

%s3.3 ###
\subsection{Proof of Theorem~\ref{thm:main-theorem}}

\label{Subsect:Proof_main_result}

In order to apply Theorem~\ref{thm:Massart}, we have to compute the
metric entropy with bracketing of each model $\mathcal{M}_{(K,\, S)}$.
This calculation is performed in the following result for which we
provide the proof in
Appendix~\ref{sec:entropy}.
\begin{prop}[Bracketing entropy of
a model] \label{prop:entropy} Let $\eta_L:\R_+ \rightarrow\R_+$
be the increasing convex function defined by
\begin{align*}
\text{Case~\ref{config:simple-multinomial}: } \eta_L(\varepsilon) &=
(1+\varepsilon)^{L+1} - 1, \\ %\label{eq:def-eta-1}
\text{Case~\ref{config:double-multinomial}: } \eta_L(\varepsilon) &=
(1+\varepsilon)^{2L+1} - 1. %\label{eq:def-eta-2}
\end{align*}

%For any choice of $K$ and $S$, $\mathcal{M}_{(K,\, S)}$ fulfills
%\ref{cond:M}.
For any $\varepsilon\in(0,1)$,
\[
H_{[\cdot]}\left(\eta_L(\varepsilon), \mathcal{M}_{(K,\, S)}, \h
\right)
\leq D_{(K,\, S)} \ln\left(\frac{1}{\varepsilon}\right) + C_{(K,\, S)},
\]
where
%
%e10 ###
\begin{equation}
\begin{split} \label{eq:def-RKS}
C_{(K,\, S)} &= \frac{1}{2} \bigg( \ln(2\pi e) D_{(K,\, S)} + \ln
(4\pi
e) \left(\1_{K\geq2} + L + (K-1) |S|\right) \\
&\quad+ \1_{K\geq2} \ln(K+1) + \sum_{l=1}^{L} \ln(A_l+1) + (K-1)
\sum
_{l\in S} \ln(A_l+1).\bigg)
\end{split}
\end{equation}
\end{prop}
The technical quantity $C_{(K,\, S)}$ measures the complexity
of a model $\mathcal{M}_{(K,\, S)}$.

The next step establishes an expression for $\phi_m$. Proof
for all subsequent results is provided in Appendix~\ref{sec:phi}.
\begin{prop} \label{prop:phi-m} For any choice of $m=(K,\, S)$,
the function $\phi_m$ defined on $(0,\eta_L(1)]$ by
\[
\phi_m(\sigma) = \left(2\sqrt{\ln2}\,\sqrt{D_{(K,\, S)}} + \sqrt
{C_{(K,\, S)} - D_{(K,\, S)} \ln\eta_L^{-1}(\sigma)} \right)\,
\sigma
\]
fulfills~\ref{cond:i} for $\sigma\leq\eta_L(1)$.
\end{prop}
To avoid more complicated expressions, we do not define $\phi_m$ for
$\sigma$ bigger than $\eta_L(1)$. A condition on $\xi$ therefore
appears in the following lemma:
\begin{lem} \label{lem:sigma-less-than-eta} For both Case \ref
{config:simple-multinomial} and Case~\ref{config:double-multinomial},
for all $n \geq1$, if $\xi =\frac{4\sqrt{LA_{\max}}}{\eta_L(1)}\leq1$
the solution $\sigma _m$ of (\ref{eq:def-sigma-m}) satisfies
$\sigma_m<\eta_L(1)$.
\end{lem}
The condition appearing in Lemma~\ref{lem:sigma-less-than-eta}
is fulfilled unless $L$ is very small, which is not the case for
the usual applications.

We can deduce an upper bound for $\sigma_m$ based on Proposition~\ref
{prop:phi-m} with a similar reasoning to \cite{MaugisMichel2008thm}.
First, $\sigma _m\leq\eta_L(1)$ implies $\eta_L^{-1}(\sigma_m)\leq1$,
and we obtain the lower bound $\sigma_m \geq\widetilde{\sigma}_{m}$,
where
%
%e11 ###
\begin{equation} \label{eq:sigma-tilde}
\widetilde{\sigma}_{m} = \frac{1}{\sqrt{n}} \left(2\sqrt{\ln2}\,
\sqrt
{D_{m}} + \sqrt{C_{m}} \right).
\end{equation}
This can be used to get an upper bound
%
%e12 ###
\begin{equation} \label{eq:sigma-upper}
\sigma_m \leq\frac{1}{\sqrt{n}} \left(2\sqrt{\ln2}\,\sqrt
{D_{m}} +
\sqrt{C_{m} - D_{m} \ln\eta_L^{-1}\left(\widetilde{\sigma
}_{m}\right)}
\right).
\end{equation}

We then choose the weights $x_m$. For values bigger than $n \sigma_m^2$,
this will change the shape of the penalty in Theorem~\ref{thm:Massart}.
We define
\begin{equation*} \label{eq:x-m}
x_m = (\ln2) D_{m}.
\end{equation*}
The following lemma shows that this is a suitable choice.
\begin{lem} \label{lem:sum-x-m}
For any model $\mathcal{M}_{m}$, with $m\in\colmodels$ as above, let
$x_m = (\ln2) D_{m}$. Then
\[
\sum_{m\in\colmodels} e^{-x_m} \leq(3/4)^L.
\]
\end{lem}

We must lower bound $\eta_L^{-1}(\widetilde{\sigma}_{m})$ to express
the penalty function, which is accomplished in the following lemma.
\begin{lem} \label{lem:pen}
Using the preceding notations,
\[
\sigma_m^2 + \frac{x_m}{n} \leq\Bigg(5+\sqrt{\max\left(\frac
{\ln n+\ln
L}{2},\; \frac{\ln2}{2}+\ln L\right)}\Bigg)^2 \frac{D_{(K,\, S)}}{n}.
\]
\end{lem}

We finally use Theorem~\ref{thm:Massart} to complete the proof of
Theorem~\ref{thm:main-theorem}.

%s4 ###
\section{Practical application}

\label{Sect:In_practice}

% Il faut aussi parler de la possibilit\'e de prendre en compte les
%nouveau modeles dans l'implementation

In real datasets, the number $A_l$ of all possible modalities for each
variable $X^l$ is not necessarily known. However, the observed number
can be used instead. In fact, the MLE estimator
selects a density with null weight on non-observed alleles. Then,
in each model $\mathcal{M}_{(K,\, S)}$, an approximated ML-estimator
can be computed thanks to the Expectation-Maximization (EM) algorithm of
\cite{Dempster1977}.

We use the same EM strategy as \cite{ToussileGassiat2009} to
avoid a local maximization of the likelihood: We run a certain number
($15$ by default) of iterations in the EM algorithm from several
($10$ by default) randomly chosen parameter points and perform a long
EM run of the best candidate in terms of likelihood.

% sets of frequencies were simulated uniformly at random on
% the corresponding simplex; they were used to initialize several EM
% algorithms in parallel, and after a certain number of steps ($15$
% by default), only the best candidate was retained to continue the
% EM algorithm. %In our simulations with this setting we did not
%observed bad calculations of the MLE due to local maximization, while
%bad calculations were frequent when the initialization frequencies
%were taken according a Dirichlet $D(1/2, \ldots, 1/2)$ distribution.

Two other points that have to be addressed before obtaining the final
estimator $\widehat{P}_{(\widehat{K}_n,\ \widehat{S}_n)}$ concern
the choice of the penalty function and the sub-collection of models
among which to select the optimal model. These two points are discussed
in Subsections~\ref{Subsect:slope_heuristics} and \ref
{Subsect:sub-collection}.
Simulations are presented in Subsection~\ref{sect:Simulated_experiments}.

%s4.1 ###
\subsection{Slope heuristics and dimension jump}

\label{Subsect:slope_heuristics}%TODO s'inspirer de la page 126 de la
%th
%\`ese de Lebarbier

Theorem~\ref{thm:main-theorem} suggests to use a penalty function of
the shape given in equation~(\ref{eq:shape.penalty}), where modulo is
defined as a multiplicative parameter $\lambda$ that has to be
calibrated. Slope heuristics, as presented in \cite{BirgeMassart2007}
and \cite {ArlotMassart2009}, provide a practical method to find an
optimal penalty $\pen_{\mathrm {opt}}(m) = \lambda_{\mathrm{opt}}
D_m/n$. These heuristics are based on the conjecture that a minimal
penalty $\pen_{\min}(m) =\lambda_{\min} D_m/n$ exists that is required
for the model selection procedure: when the penalty is smaller than
$\pen_{\min}$, the selected model is one of the most complex models,
and the risk of the selected estimator is large. In contrast, when the
penalty is larger than $\pen_{\min}$, the selected model is
considerably less complex. Thus, the optimal penalty is close to twice
the minimal penalty:
\begin{equation*}
\pen_{\mathrm{opt}}\left(m\right)\approx2\lambda_{\min} \,\frac{D_m}{n}.
\end{equation*}

\begin{sloppypar}
An explanation of the heuristics behind this factor $2$ can be found in
\cite{MaugisMichel2008slope}, for instance. The name ``slope
heuristics'' is derived from $\lambda_{\min}$, which is the slope of
the linear regression $\gamma_n(\widehat{P}_{m})\sim D_m/n$ for a
certain sub-collection of the most competitive models $m$. For example,
in Figure~\ref{Chap3.Fig:slope.heuristics_a} below, models containing
the true one $\mathcal{M}_{(K_0,\ S_0)}$ exhibit a slope. This example
also illustrates that slope heuristics are appropriate in our modeling
context.\looseness=1
\end{sloppypar}

%f1 ###
\begin{figure}[t]
\centering
\subfigure[$\widehat{\lambda}_{\min}\approx0.88$ by linear regression
on sub-models.]{\label{Chap3.Fig:slope.heuristics_a}
\ifpdf
\includegraphics[width=0.45\textwidth]{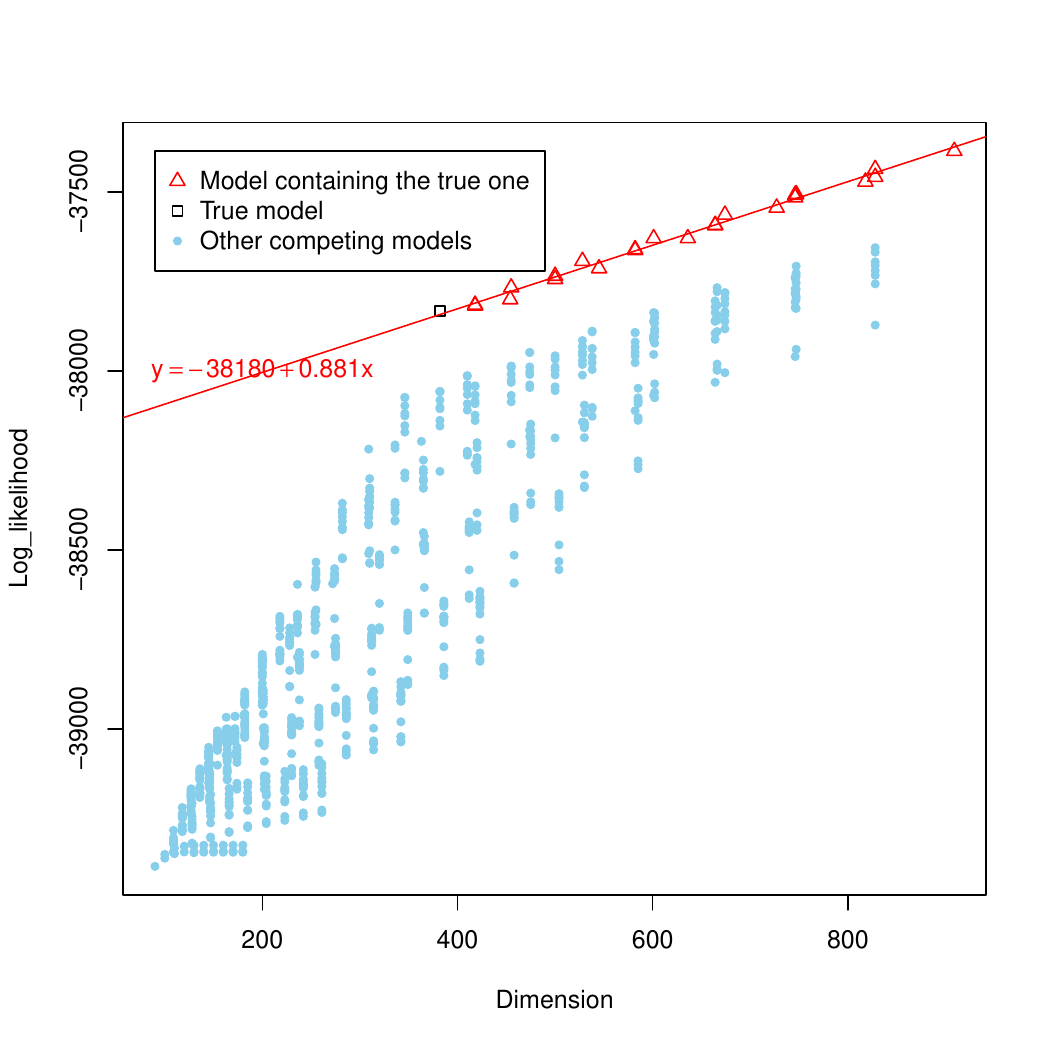}
\else
\includegraphics[width=0.45\textwidth]{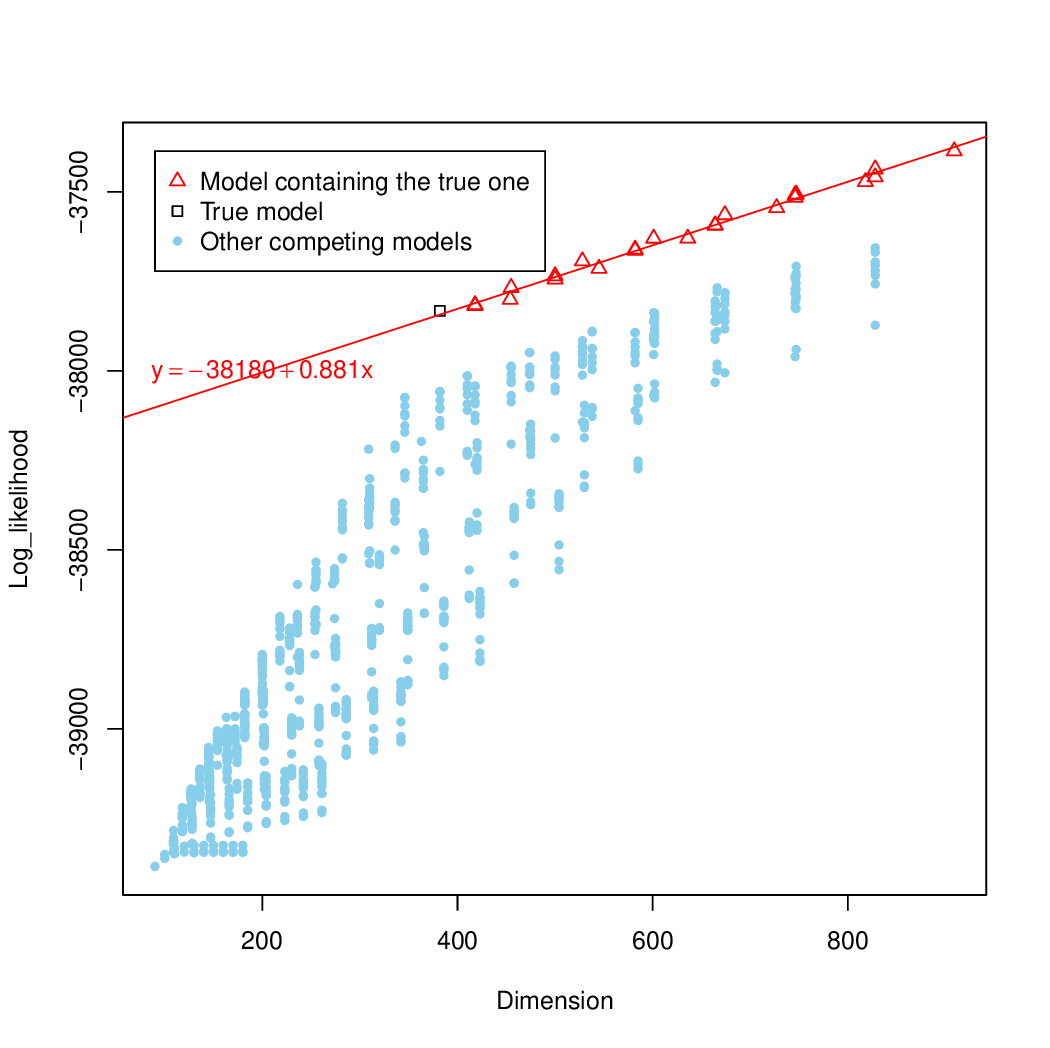}
\fi
}\quad
\subfigure[The biggest jump occurs around $5.10$, but the sliding
window selects the more reasonable value $\widehat{\lambda}_{\min
}\approx0.9$.]{\label{Chap3.Fig:slope.heuristics_b}
\ifpdf
\includegraphics[width=0.45\textwidth]{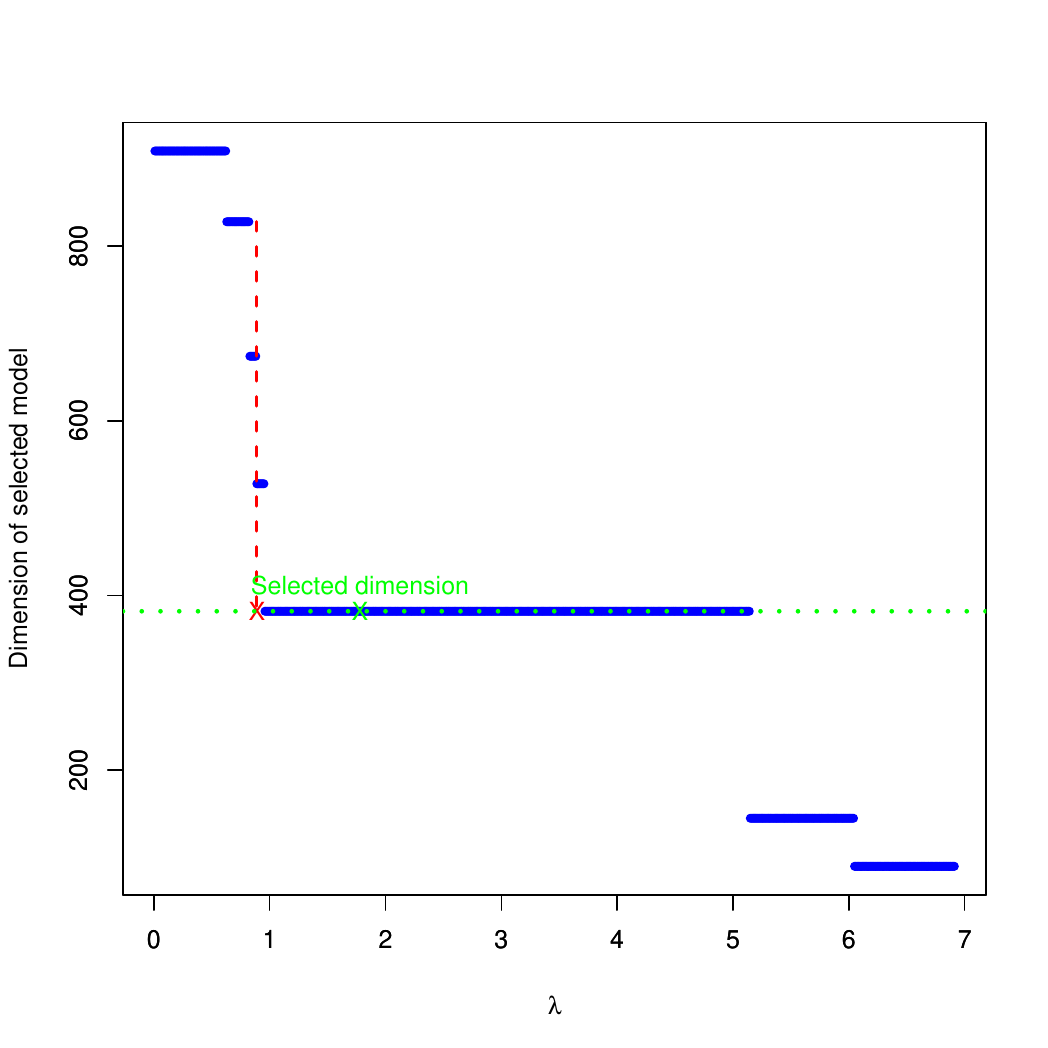}
\else
\includegraphics[width=0.45\textwidth]{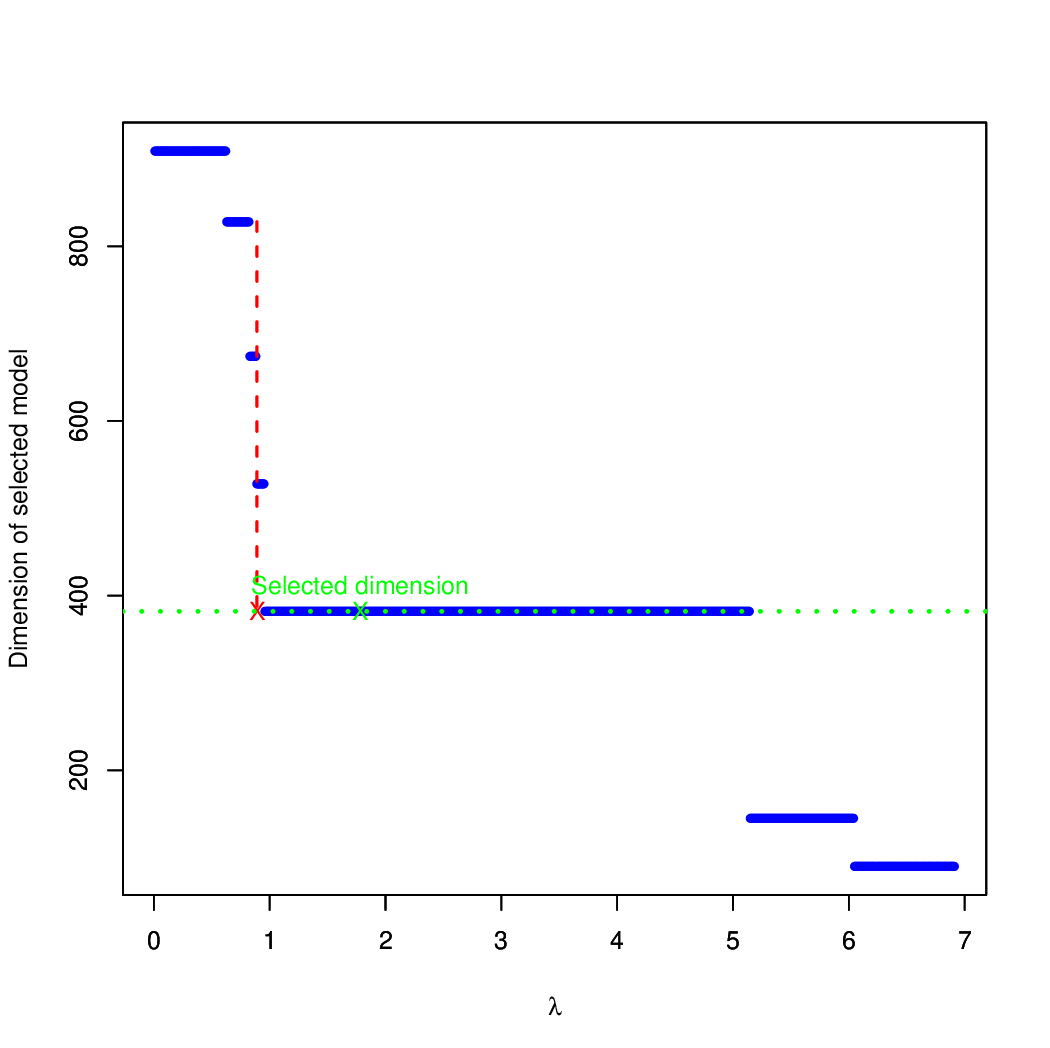}
\fi
}
\caption{Two ways to compute the slope for a simulated sample of $1000$
individuals with $8$ clustering loci among $10$ and $5$ populations.
Models are explored via the %modified
backward-stepwise method described in Subsection~\ref{Subsect:sub-collection};
the number of clusters $K$ ranges from $1$ to $10$. The size of
the sliding window is $0.10$.}
\label{Fig:slope.heuristics.1}
\end{figure}

Instead of using linear regression to estimate $\lambda_{\min}$, we use
another method (generally referred to as the dimension jump method) to
detect the biggest jump on the selected dimension with respect to the
candidate values of $\lambda$. In practice, we would assume a
reasonable grid $\lambda_1<\cdots<\lambda_{n_{\lambda}}$ of
$n_{\lambda}$ candidate estimates of $\lambda_{\min}$ and a
sub-collection $\mathcal{C}_{ex}$ comprising the most competitive
models. Each $\lambda_i$ leads to a selected model $\widehat{m}_i$ with
dimension $D_{\widehat{m}_i}$. If $D_{\widehat{m}_i}$ is plotted as a
function of $\lambda_i$, $\lambda_{min}$ is expected to lie at the
position of the biggest jump.\looseness=1

However, Fig.~\ref{Chap3.Fig:slope.heuristics_b} illustrates an
important point: in this example the biggest jump occurs at $\lambda
\approx5.1$, but the optimal value of $\lambda_{\min}$ is around $0.9$,
which corresponds to several successive jumps. %This does not necessarily
%come from bad computations of the ML-estimators (via EM); considering
%some very high dimension models could also improve the calculus of the
%slope.
We propose an improved version of the dimension jump method of \cite
{ArlotMassart2009}
based on a sliding window: on the axis of $\lambda$, we consider the
sum of all jumps in a sliding window of size $h>0$. In Algorithm~\ref
{alg:calibration} below, which describes the procedure, $n_h$ denotes
the number of candidate values of $\lambda_{\min}$ in the sliding window.
%An example of the application of this variant of dimension jump method
%is given in Figure~\ref{Chap3.Fig:slope.heuristics_c}.
We do not claim that this improves slope heuristics per se; we merely
note that the proposed procedure improves the stability of the method
in our simulations. In practice, following repeated trials, we choose a
window of size $h = 0.10$.\looseness=1

\begin{algorithm}[t]
  \caption{Penalty Calibration $\left(\mathcal{C}_{ex}, \left(\lambda_i\right)_{i=1,\ldots,n_{\lambda}}, n_h\right)$}
  \label{alg:calibration}
  \begin{algorithmic}
    \FOR{$i=1$ to $n_\lambda$}
  \STATE $\displaystyle \widehat{m}_i \leftarrow  \argmin_{m\in\mathcal{C}_{ex}}\left\{ \gamma_n\left({P}_m\right)
  +\lambda_i D_m/n\right\}$
      \ENDFOR
    \STATE $\displaystyle i_{end} \leftarrow \min \argmax_{i\in\left\{n_h+1,\ldots,n_{\lambda}\right\}}\left\{D_{\widehat{m}_{i-n_h}}-D_{\widehat{m}_{i}}\right\}$
    \STATE $i_{init} \leftarrow \max \Big\{ j \in \left[i_{end}-n_h, i_{end}-1\right], D_{\widehat{m}_{j}}-D_{\widehat{m}_{i_{end}}} = D_{\widehat{m}_{i_{end}-n_h}}-D_{\widehat{m}_{i_{end}}}\Big\}$
    \STATE $\displaystyle \widehat{\lambda}_{\min} \leftarrow \frac{\lambda_{i_{init}}+\lambda_{i_{end}}}{2}$
    \RETURN $\widehat{\lambda}_{\min}$
  \end{algorithmic}
\end{algorithm}

% \subsection{Sub-collection of models used to calibrate the penalty}
%s4.2 ###
\subsection{Sub-collection of the most competitive models}

\label{Subsect:sub-collection}

For a given maximum number of clusters $K_{\max}$, the number of
competing models is equal to $1+(K_{\max}-1)*(2^{L}-1)$. Because this
is a very large number in most situations, it would be very laborious
to consider the total number of potentially applicable models to
calibrate the parameter $\lambda$. Nevertheless, a sufficient number of
models is necessary to ensure a clear jump in the selected dimension
sequence. We therefore consider the modified backward-stepwise
algorithm proposed in \cite{ToussileGassiat2009}, which enables us to
gather the most competitive models among all possible $S$ for a given
number of clusters $K$ and a given penalty function $\pen_{n}$. This
algorithm also offers the possibility to add a complementary
exploration step based on a similarly modified forward strategy: we
refer to this algorithm as $explorer(K,\ \pen_{n})$.

Because the final penalty during the exploration step is unknown, we
consider a reasonable grid $\frac{1}{2}=\lambda_1<\cdots<\lambda
_{n_{\lambda}}=\ln n$ containing both penalty functions associated with
$\AIC$ and $\BIC$. Each value $\lambda_i$ is associated with a penalty
function $\pen _{\lambda_i}$. We launch $explorer(K,\
\pen_{\lambda_i})$ for all values of $K$ in $\{ 1,\ldots,\ K_{\max}\} $
and for all values of $\lambda_i$ of the grid; we then gather the
explored models in $\mathcal{C}_{ex}$. This sub-collection appears to
contain the most competitive models and it was therefore used to
calibrate $\lambda$.

%s5 ###
\section{Simulations}

\label{sect:Simulated_experiments}

% The procedure we propose in this paper uses a data-driven calibration
% of the penalty function.
Our proposed procedure is implemented in the software
\texttt{MixMoGenD} (Mixture Model for Genotypic Data), which already
offers a selection procedure based on the asymptotic criteria $\BIC$
and $\AIC$ \citep{ToussileGassiat2009}. Numerical experiments
with simulated datasets are performed to assess the performance of the
new non-asymptotic
criterion with respect to $\BIC$, $\AIC$, and the Integrated Completed
Likelihood ($\ICL$) \citep{biernacki2000}. %We used the $\BIC$ type
%approximation of $\ICL$ that leads to $\ICL\left(K,S\right)=\BIC
%\left(K,S\right)+\ENT\left(\widehat{P}_{(K,S)}\right)$, where $\ENT$
%is an entropy term given by $\ENT\left(\widehat{P}_{(K,S)}\right)=-
%\sum_{i=}^n\sum_{k=1}^K\widehat{z}_{ik}\ln\left(\widehat{t}_{ik}
%\right)$. %The penalty functions of these last criteria are
%respectively defined by
%\begin{align*}
% \pen_{\BIC}\left(m\right)&=\frac{\ln n}{2n}\, D_m;\\
% \pen_{\AIC}\left(m\right)&=\frac{1}{n}\, D_m.
%\end{align*}

We set up two series of experiments to simulate multilocus genotypic
data from diploid organisms (Case~\ref{config:double-multinomial}).
%For such data, each variable corresponds to a set of $2$ unordered
%nominal variables called alleles.
Consistency behaviors of the competing criteria are then evaluated
based on the first series, which examines how the main features of the
true model are retrieved as the sample size increases. The second
series compares the risks of the selected estimators from an oracle
perspective.

%s5.1 ###
\subsection{Consistency behaviors}

\label{sect:consistency}

\begin{figure}[t]
\ifpdf
\includegraphics[width=0.49\textwidth]{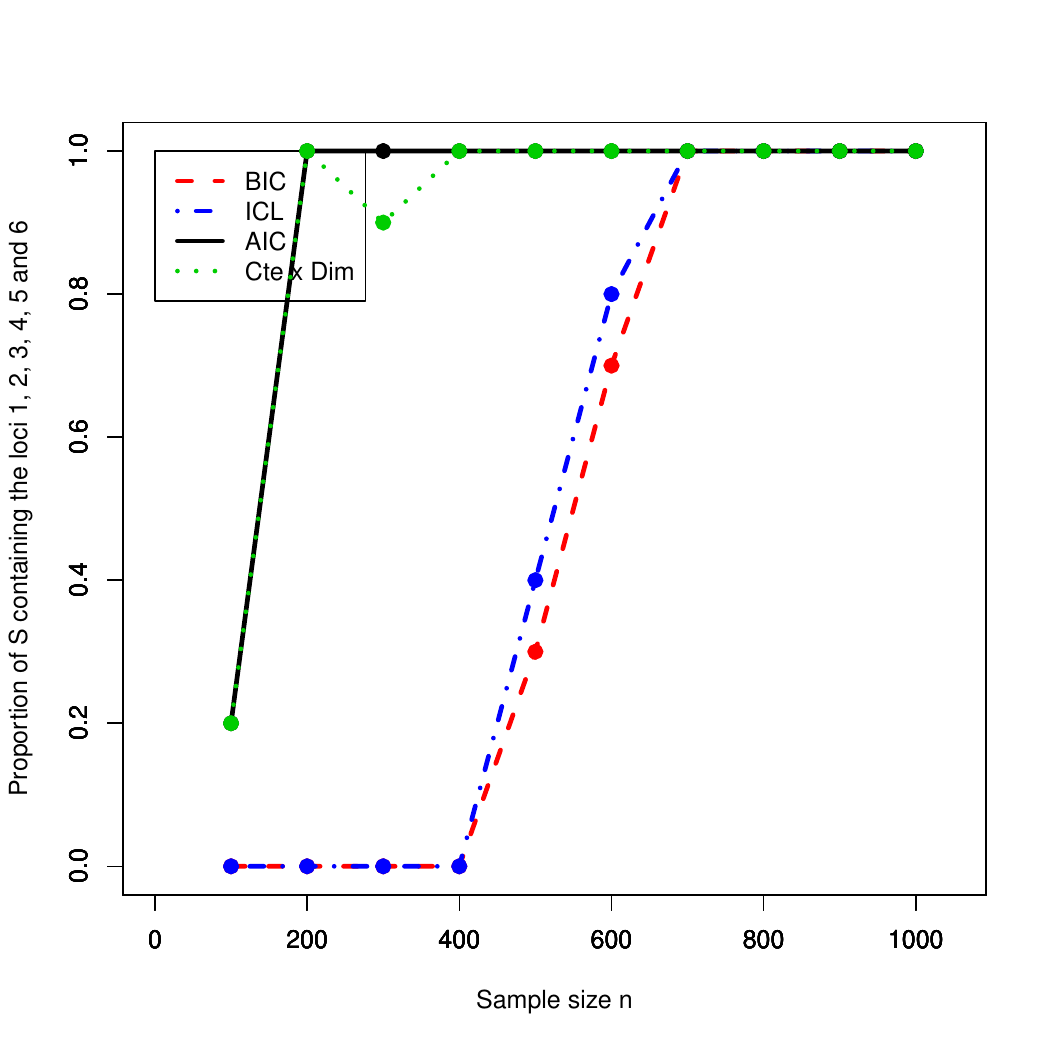}
\includegraphics[width=0.49\textwidth]{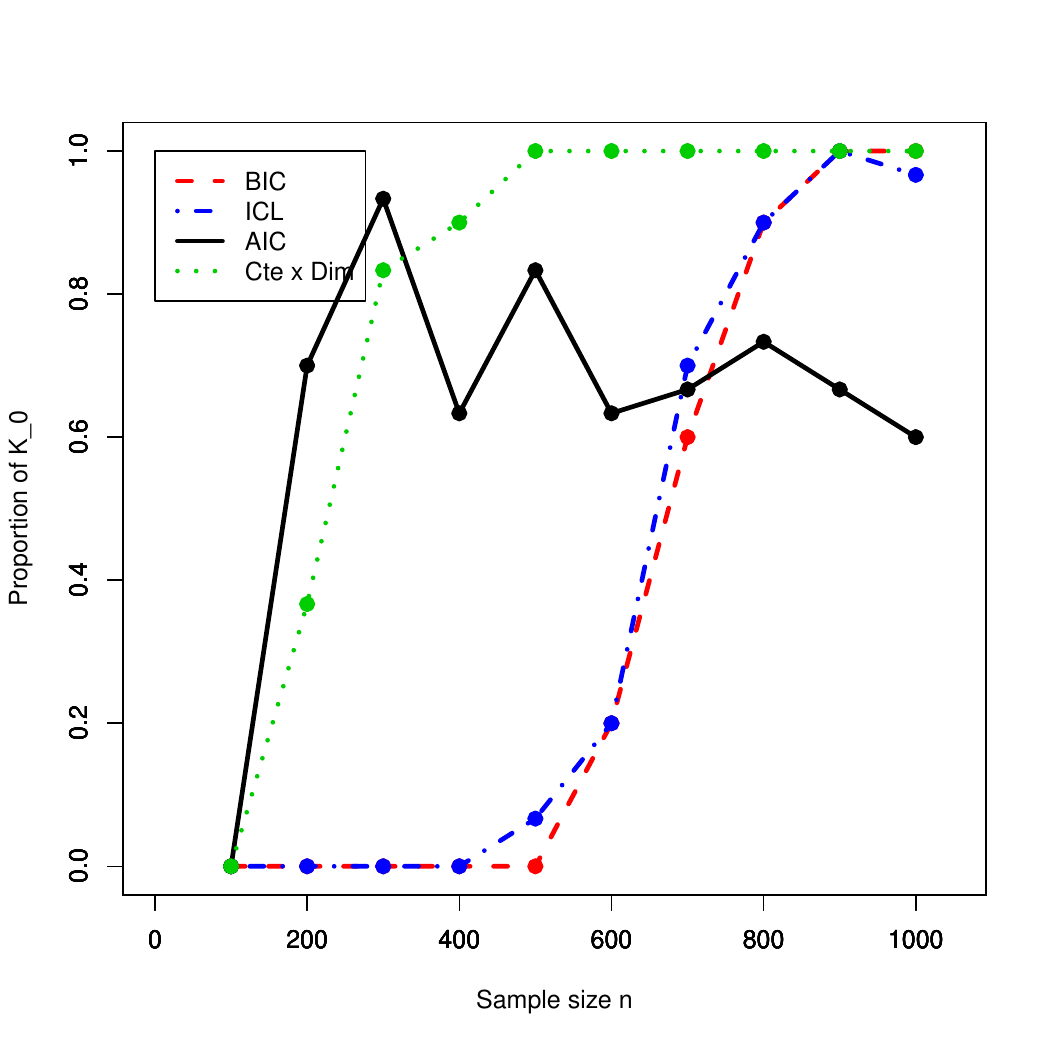}
\else
\includegraphics[width=0.49\textwidth]{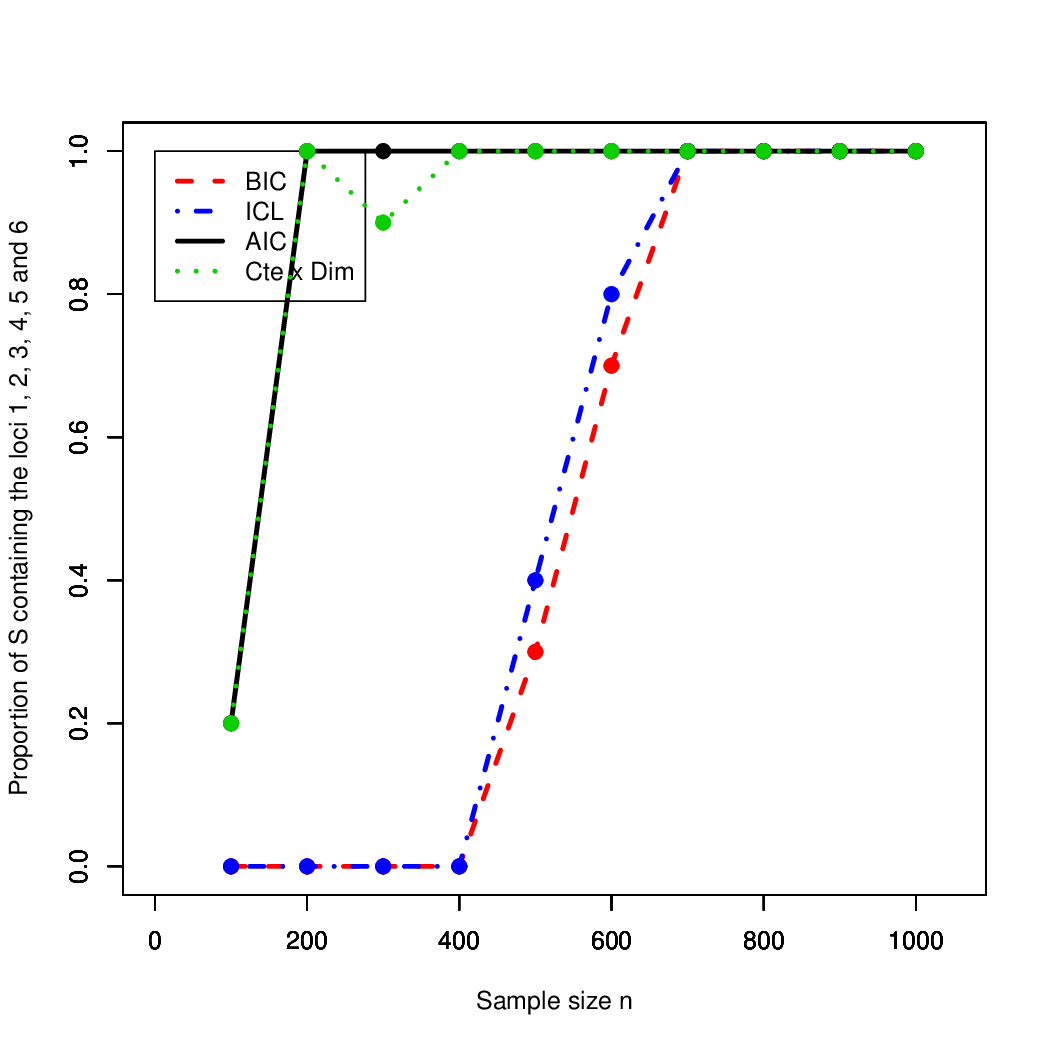}
\includegraphics[width=0.49\textwidth]{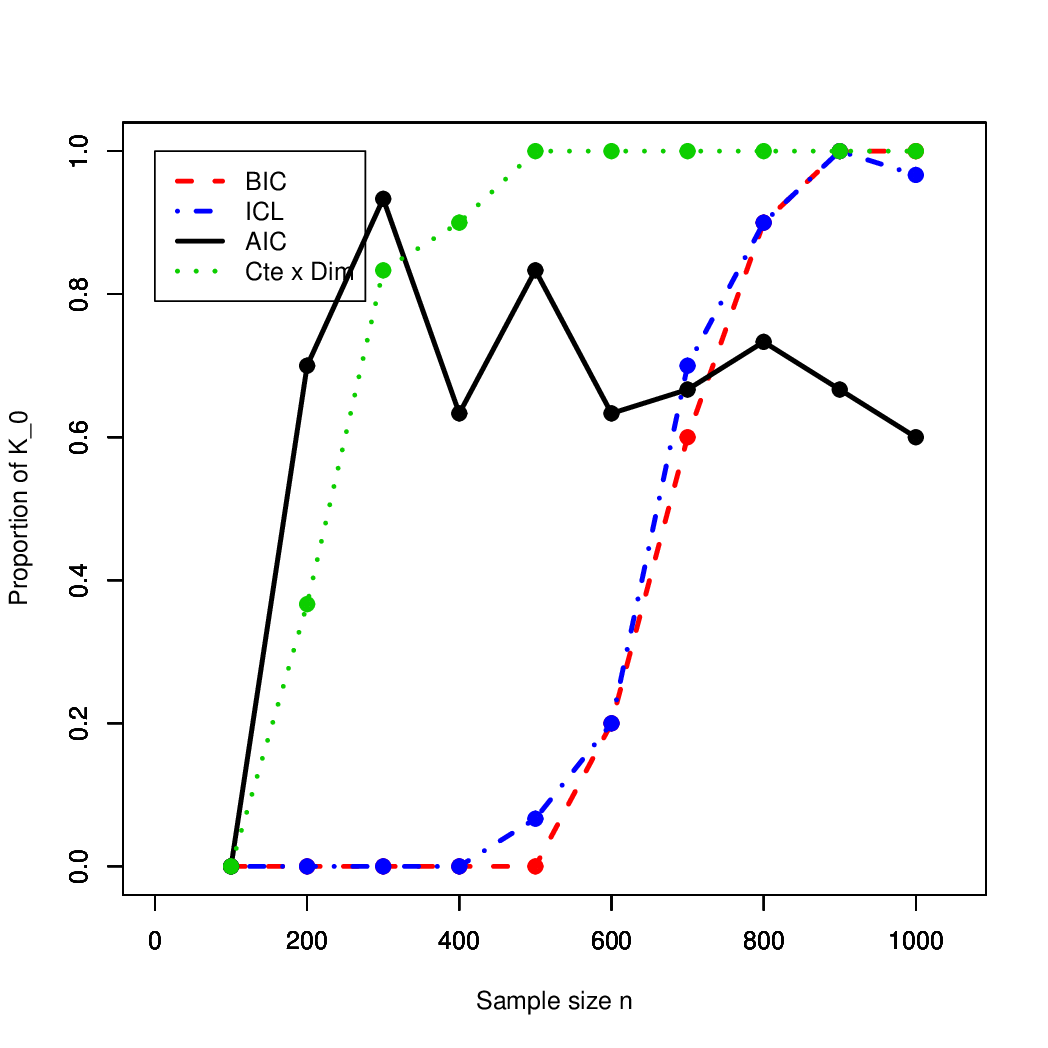}
\fi
\caption{Left: Graph showing the proportion
of selected models with $\widehat{S}_n$" containing ${1,\ldots,6}$ with
respect to sample size $n$. Right: Graph showing the proportion of
selected models with $\widehat{K}_n=K_0$ with respect to sample size $n$.}
\label{Fig:S_vs_n-K_vs_n}
\end{figure}

We consider a setting of $L=10$ variables with $10$ categories each.
Each dataset is simulated as a mixture of $K_0=5$ populations in equal
proportions. The simulation parameters are chosen so that the
differentiation between populations, as measured by a population
genetics parameter $F_{st}$ (a measure of genetic differentiation),
decreases with the variable rank. Populations are distinctly separated
for the first $6$ variables; for the next $2$ variables populations are
poorly differentiated; the last $2$ variables follow the uniform
distribution for all populations. The complete parameter is available
at \url{http://www.math.u-psud.fr/\textasciitilde toussile/}. The
overall differentiation occurs in a range considered difficult for
clustering of such data \citep{Latch2006}. We examine different values
of the sample size $n$ in $[100,\ 1000]$, and $30$ datasets are
simulated for each value. Results are summarized in
Fig.~\ref{Fig:S_vs_n-K_vs_n}.
% and Table~\ref{Tab:Select.K}.
%The left panel gives the proportion of selecting the subset $
%\widehat{S}_n$ of clustering variables containing the first $6$
%variables, which are the most genetically differentiated variables.
%The right panel gives the proportion of selected models with $
%\widehat{K}_n=K_0$. %(rÂŽpÂŽtition avec la lÂŽgende).

%f2 ###

%f3 ###
\begin{figure}[b!]
\begin{center}
\ifpdf
\includegraphics[width=0.6\textwidth]{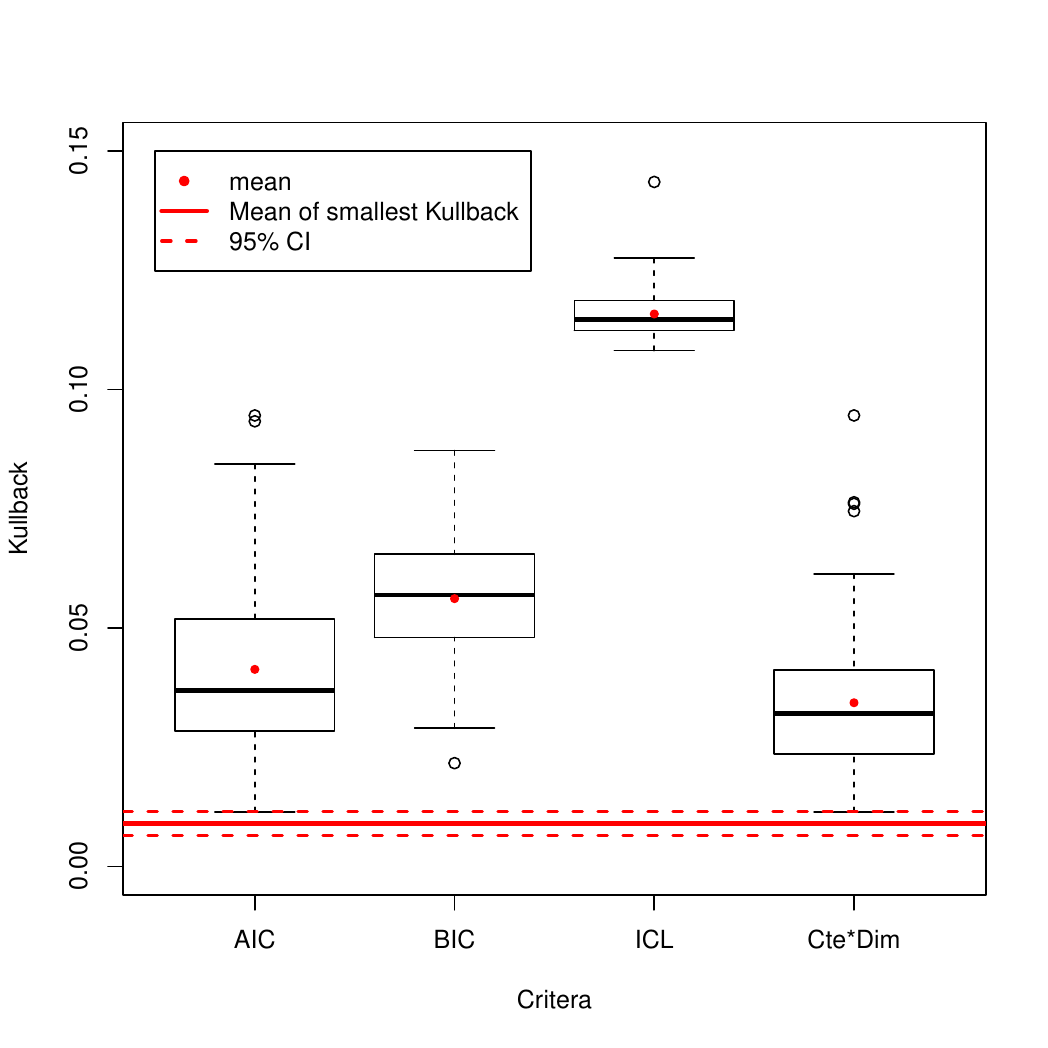}
\else
\includegraphics[width=0.6\textwidth]{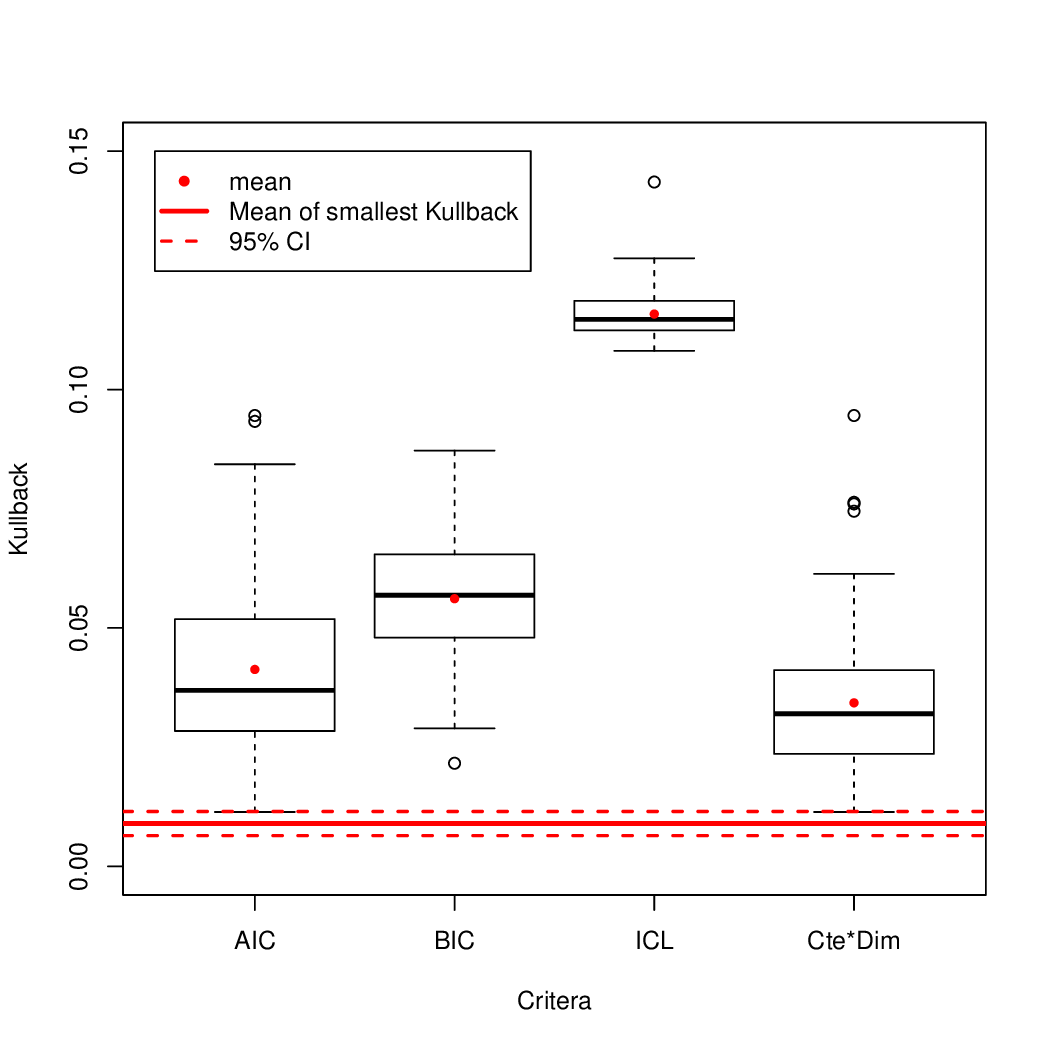}
\fi
\caption{Box plots of the Kullback risk (estimated with a Monte Carlo
procedure)
of the selected estimator. The red line corresponds to the mean (and a
$95\%$ confidence interval) of the smallest risk obtained on competing
estimators from each dataset. $\CteDim$ denotes the new criterion
with a data-driven calibration of the penalty function.}
\label{Fig:Percentage.oracle}
\end{center}
\end{figure}

We observe similar behaviors of $\BIC$ and
$\ICL$ in these experiments: both criteria perform poorly for the
selection of variables and the classification of small sample sizes. In
fact, \cite{Nadif1998}
have pointed out that $\BIC$ requires a large sample size to
reach its asymptotic behavior in a discrete framework. The high
variability of the dimensions
of the competing models, which cancels the contribution of the entropy
term in $\ICL$, may explain the similar behavior of $\BIC$ and $\ICL$.
In contrast, $\AIC$ and the newly proposed criterion are most suited
to the selection of variables for both small and large sample sizes.
The new criterion also performs well for the selection of the number
of components for both small and large sample sizes, but $\AIC$ overestimates
the number of components for large sample sizes (from $n=400$). %, see
%Table~\ref{Tab:Select.K}).
As expected, the data-driven calibration of the penalty function globally
improves the performance of the selection procedure and consequently
provides an answer to the question ``Which penalty for which sample
size?''

Small variations in the results obtained for small sample sizes may
occur from one run to another. In fact, the EM algorithm may fail to
identify the global maximum for such sample sizes, in particular for
models of larger dimensions. This is probably the case for some
datasets of size $n\leq300$; the number of free parameters in our
simulated model is $\geq310$.\vspace*{-3pt}

%s5.2 ###
\subsection{Oracle performance}

\label{sect:oracle_experiment}

As previously mentioned, the new criterion is designed in a density
estimation framework. The following section compares the risks of the
selected estimators. Our simulations consist of $101$ datasets with
$L=6$ variables, $3$ categories for each variable, and $K_0=3$
components in equal proportions. The simulation parameters are chosen
in such a way that the differentiation between the components is
significant for the first $3$ variables and very small for the
$4^\text{th}$ and $5^\text{th}$ variables, whereas the $6^\text{th}$
variable follows the uniform distribution for all components. Thus, the
true model is defined by $K_0=3$ and $S_0=\{ 1,\ 2,\ 3,\ 4,\ 5\} $. The
complete parameter is available at
\url{http://www.math.u-psud.fr/\textasciitilde toussile/}.

The Kullback risk is estimated using a Monte Carlo procedure for $100$
simulated datasets, each with a sample size of $600$. Our results are
summarized in Fig.~\ref{Fig:Percentage.oracle} and
Table~\ref{Tab:pairwise.t.test}.\vadjust{\goodbreak}

%t1 ###
\begin{table}[t]
\caption{The p-values of one-sided pairwise Wilcoxon tests comparing
the Kullback risks:\break on average, the smallest Kullback divergence
is linked to the new criterion, followed by the divergence
of~$\AIC$}\label{Tab:pairwise.t.test}
\begin{tabular}{rrrr}
\hline
& ICL & BIC & AIC \\
\hline
BIC & 2.2e-16 & & \\
AIC & 2.2e-16 & 4.4e-10 & \\
CteDim & 2.2e-16 & 2.2e-16 & 0.0031 \\
\hline
\end{tabular}
\end{table}

Unsurprisingly, concerning the Kullback risk, the least favorable
behavior originates from $\ICL$, followed by $\BIC$. In fact, these
criteria are not designed to retrieve the minimal risk estimator. In
addition, $\ICL$ and $\BIC$ are based on asymptotic approximations and
may require large sample sizes. In contrast, the new criterion with a
data-driven calibration of the penalty function performes significantly
better (see Table~\ref{Tab:pairwise.t.test}). As stated previously,
both $\AIC$ and the new criterion are designed to find the minimizer of
the Kullback risk. Yet, similarly to $\BIC$ and $\ICL$, $\AIC$ is based
on asymptotic approximations. The new criterion is designed from a
non-asymptotic perspective, which may explain its advantage over
$\AIC$.

%s6 ###
\section{Application to real data sets}

\label{sect:real}

%s6.1 ###
\subsection{U.S. Congress voting data}

\begin{sloppypar}
The data set entitled ``1984 United States Congressional Voting Records
Database'' includes votes of the U.S. House of Representatives
Congressmen on $16$ key issues (disability, religion, immigration,
army, education, \ldots) identified by the Congressional Quarterly
Almanac (CQA) in \cite{Asuncion+Newman:2007}. This data set has $n=435$
instances ($267$ Democrats and $168$ Republicans). For each vote, three
possible responses are taken into account: for, against, and
abstention. The model selection procedure with calibration of the
penalty function is applied to these data. The maximum number of
clusters is set to $K_{\max}=10$. The selected number of clusters is
$\widehat{K}_n=6$, and the selected subset of relevant variables does
not include votes on disability and army issues. The confusion matrix
comparing the obtained partition and the Democrat/Republican
bi-partition is given in Table~\ref{Tab:RepublicansDemocrats}. More
than $91\%$ of clusters $3$ and $4$ are Republicans, whereas more than
$94\%$ of clusters $1$, $5$, and $6$ are Democrats. Republicans and
Democrats are equally represented in cluster~$2$. The subdivision of
the two main parties into various tendencies, a common occurrence in
politics, is reflected in these results.
\end{sloppypar}

%t2 ###
\begin{table}[t]
\caption{Confusion matrix comparing the obtained clusters with the
Democrat/Republican bi-partition. More than $91\%$ of clusters $3$ and
$4$ are Republicans, whereas more than $94\%$ of clusters $1$, $5$ and
$6$ are Democrats. Cl = Cluster}\label{Tab:RepublicansDemocrats}
\begin{tabular}{lcccccc}
& Cl 1 & Cl 2 & Cl 3 & Cl 4 & Cl 5 & Cl 6 \\
\hline
Republicans & $5$ & $8$ & $41$ & $111$ & $0$ & $3$\\
Democrats & $86$ & $8$ & $4$ & $7$ & $117$ & $45$\\
\hline
\end{tabular}
\end{table}

%f4 ###
\begin{figure}[t]
\subfigure[Selected dimension versus candidate constants from the
voting data: the selected constant is $\widehat{\lambda} = 3.04$,
leading to an optimal penalty $pen_{opt} = 6.08 * Dimension$.]{
\ifpdf
\includegraphics[width=2.00in,height=2.00in]{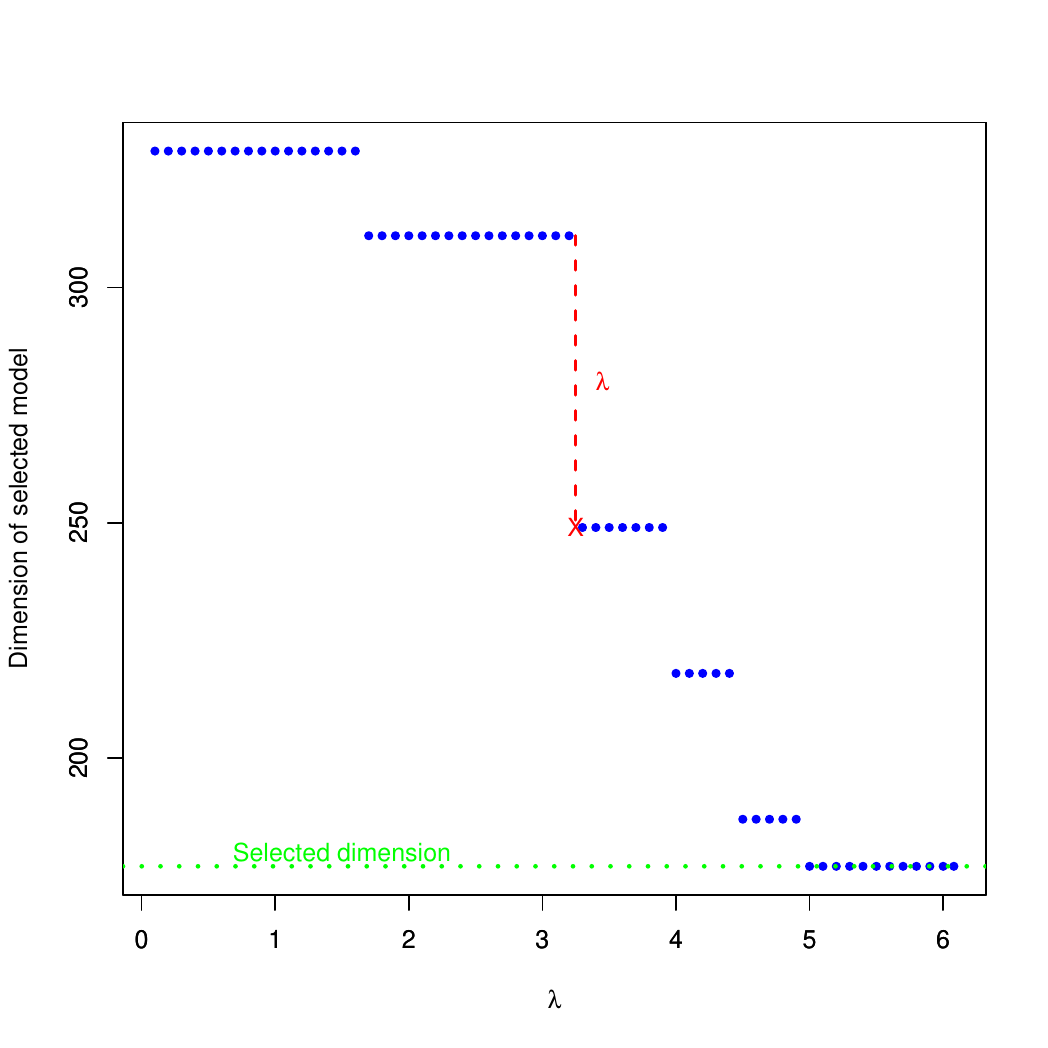}
\else
\includegraphics[width=2.00in,height=2.00in]{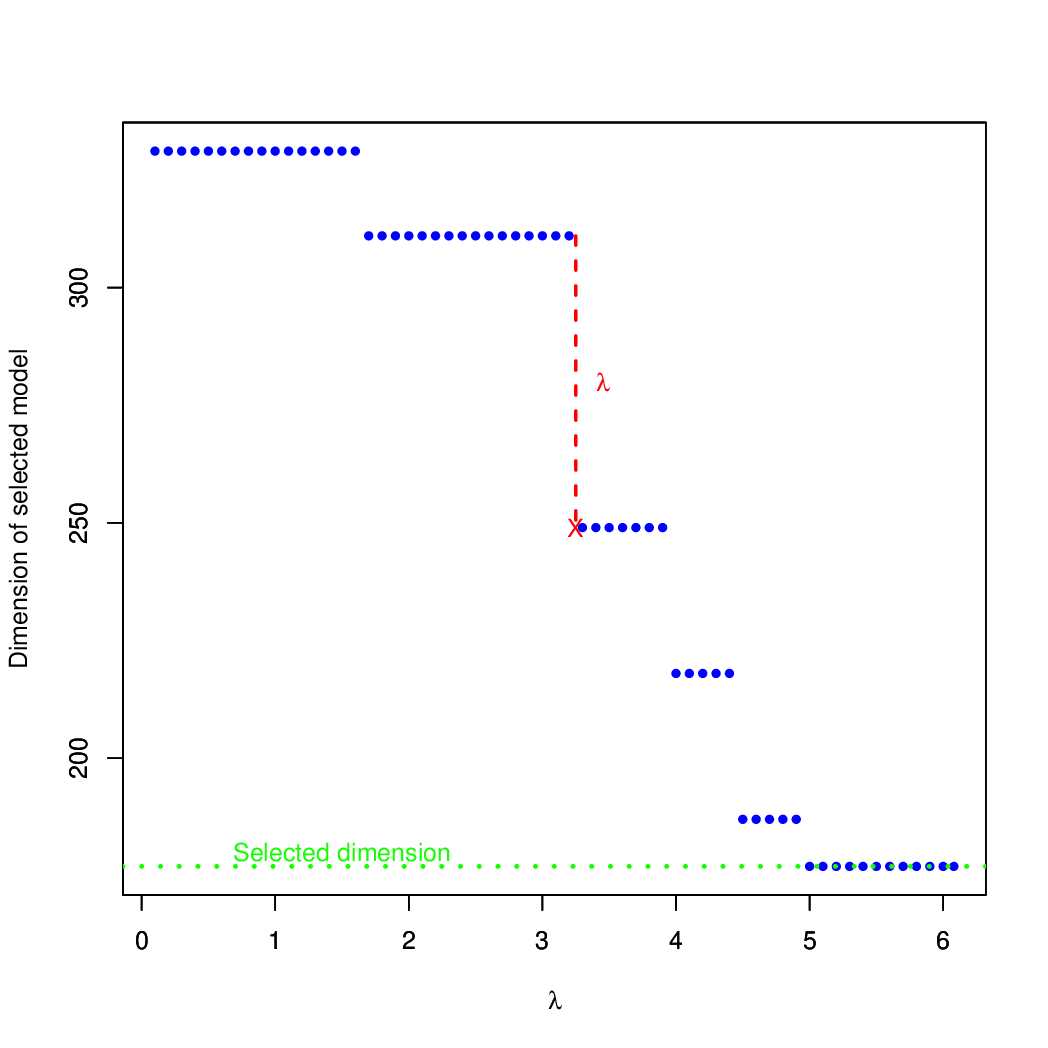}
\fi
}
\qquad
\subfigure[Log-likelihood versus Dimension of the most competitive
models from the voting data: the red line corresponds to the equation
$y=\widehat{\lambda}x + \beta$.]{
\ifpdf
\includegraphics[width=2.00in,height=2.00in]{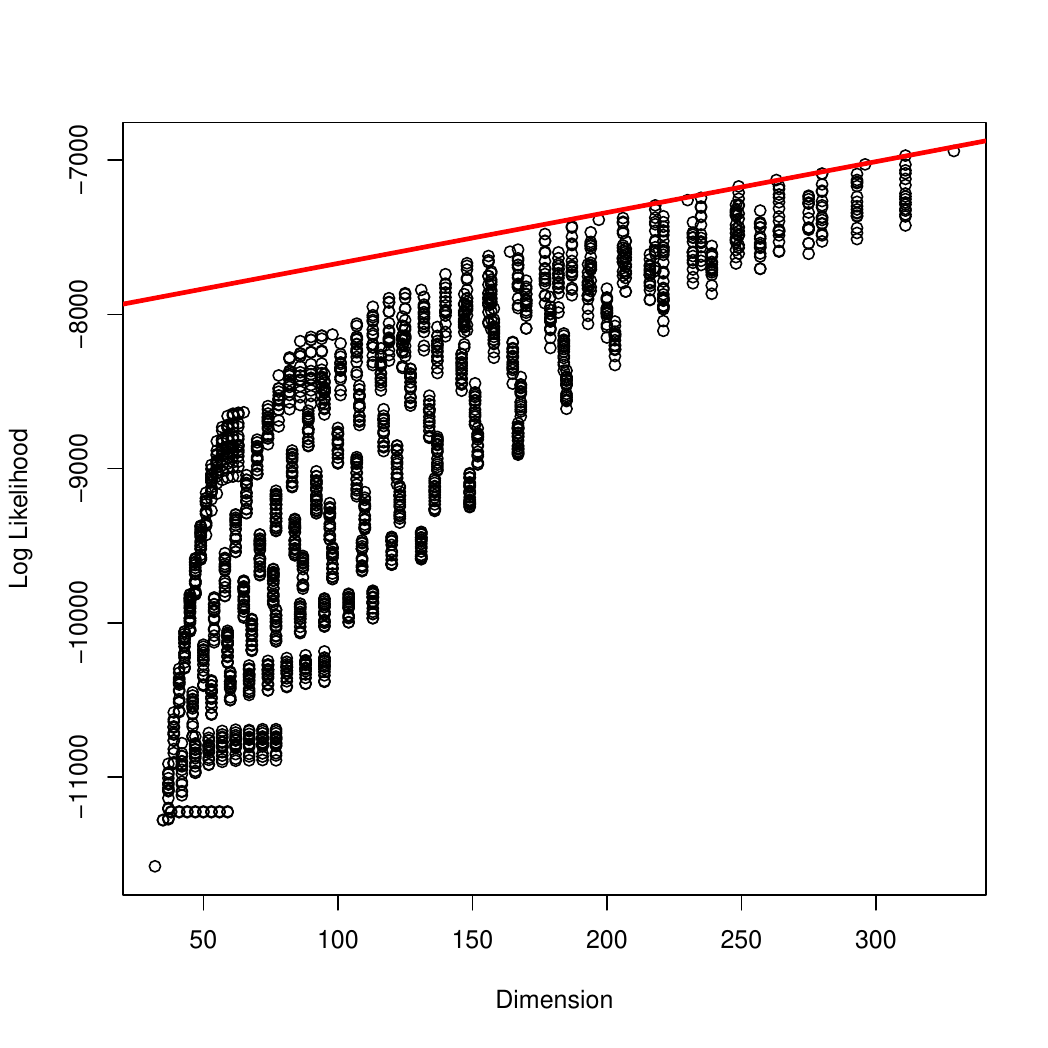}
\else
\includegraphics[width=2.00in,height=2.00in]{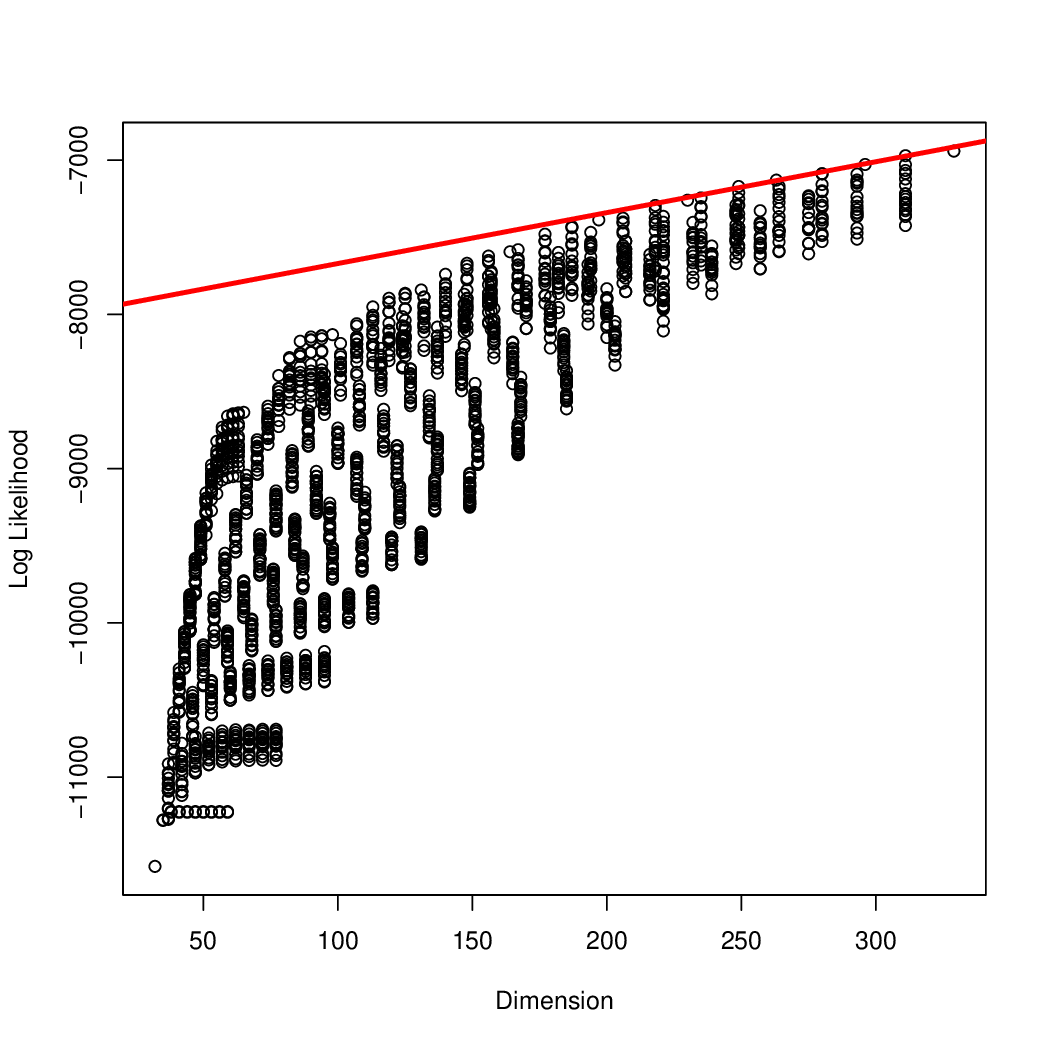}
\fi
}
\caption{Summary of the penalty calibration for voting data model selection.}
\label{Fig:ModelselectionVotingData}
\end{figure}

%s6.2 ###
\subsection{Breeds of chicken}

We consider a collection of $27$-locus genotypes from $600$ individuals
representing $20$ chicken breeds ($30$ individuals per breed).
These data have been described in \cite{Rosenberg2001a} in the context
of a clustering method evaluation of multilocus genotypes. Of the
$27$ loci, we consider $15$ that have no missing data. The data
illustrate a very common difficulty with biologic datasets: the dimensions
of the considered models are very large with respect to the number
of individuals
(note, however, that the dimensions of the competing models are large
because we have 15 variables resulting in 600x2x15 = 18 000 individual
measures). %It is difficult to validate theoretically a statistical
%study in such a case.
Nevertheless, %For these data, some competing models had a size larger
%than the size of the sample. Even in a such critical situation,
our procedure resulted in an interesting classification: $17$
clusters correspond mostly to the initial breeds, and three of the
clusters contain $2$ breeds each. The similarity between the obtained
classification and the breeds as measured by the Rand index is greater
than $98\%$. All loci are selected to be useful for clustering purposes.
In \cite{Rosenberg2001a}, the authors found $18$ clusters by using
the available $27$ variables. Their algorithm requires the user to
perform several steps. The clusters they found also corresponded mostly
to the initial breeds.

%f5 ###
\begin{figure}[htb]
\subfigure[Selected dimension versus candidate constants: th e
selected constant is $\widehat{\lambda} = 0.25$, leading to an optimal
penalty $pen_{opt} = 0.50 * Dimension$.]{
\ifpdf
\includegraphics[width=2.00in,height=2.00in]{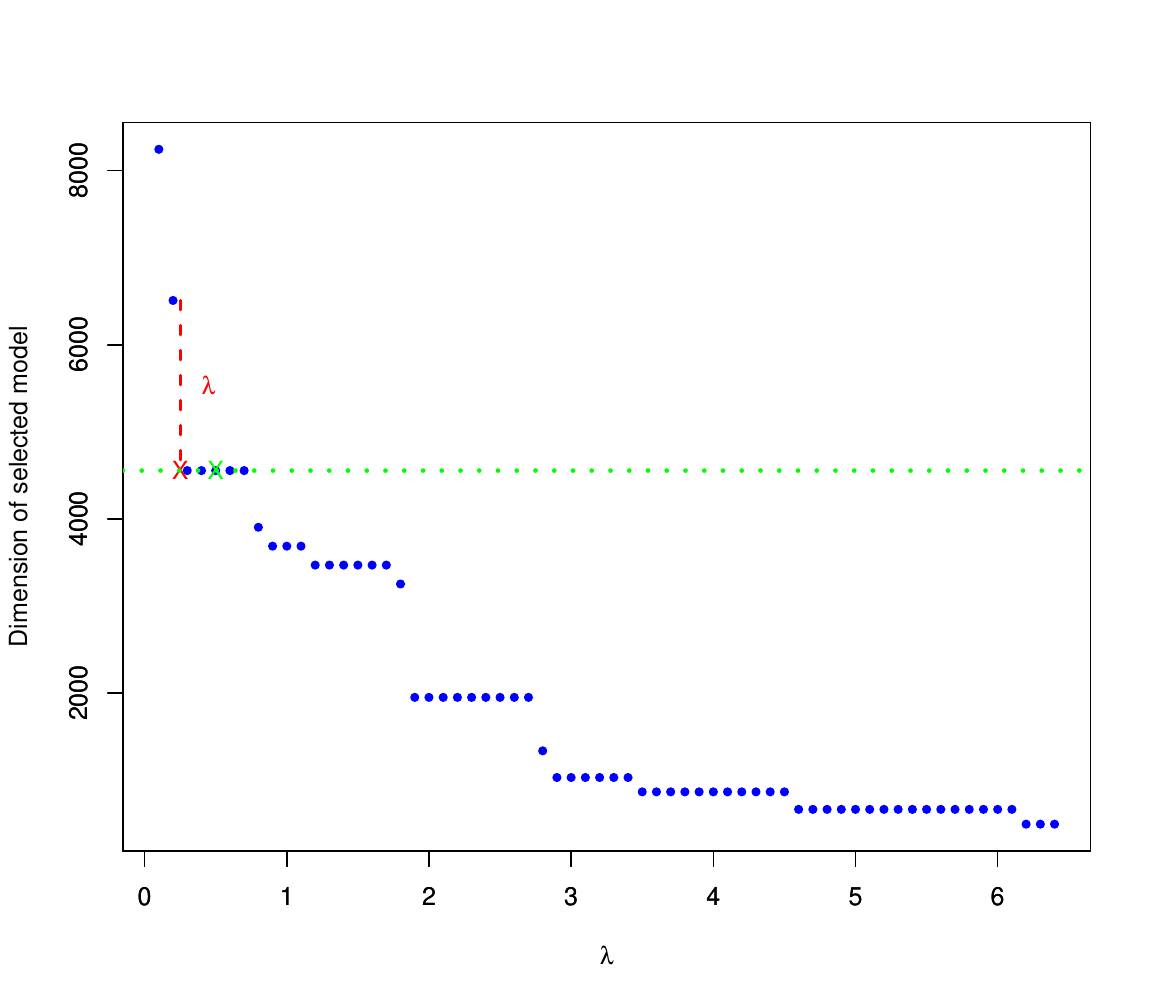}
\else
\includegraphics[width=2.00in,height=2.00in]{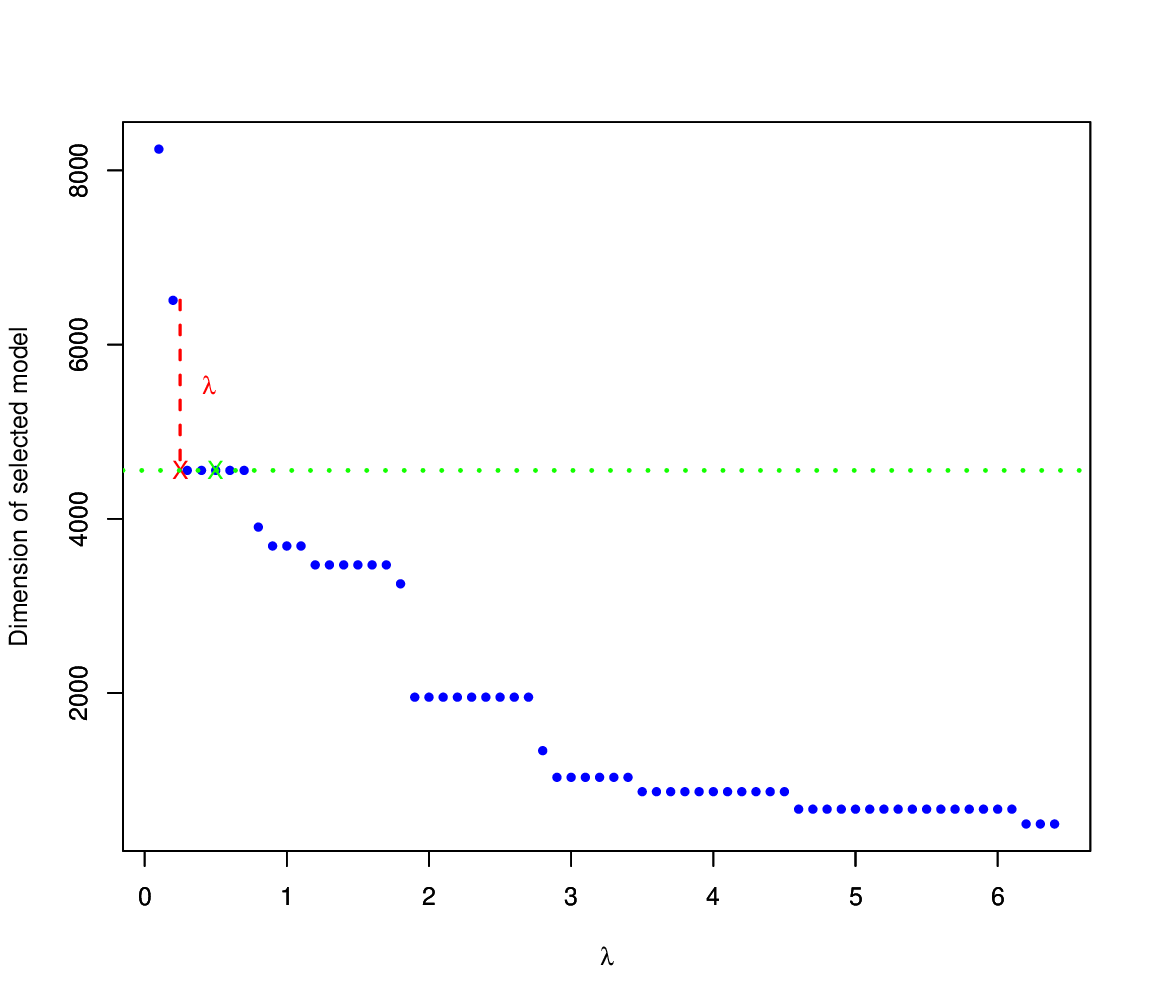}
\fi
}
\qquad
\subfigure[Log-likelihood versus Dimension of the most competitive
models from the voting data: the red line corresponds to the equation
$y=\widehat{\lambda}x + \beta$.]{
\ifpdf
\includegraphics[width=2.00in,height=2.00in]{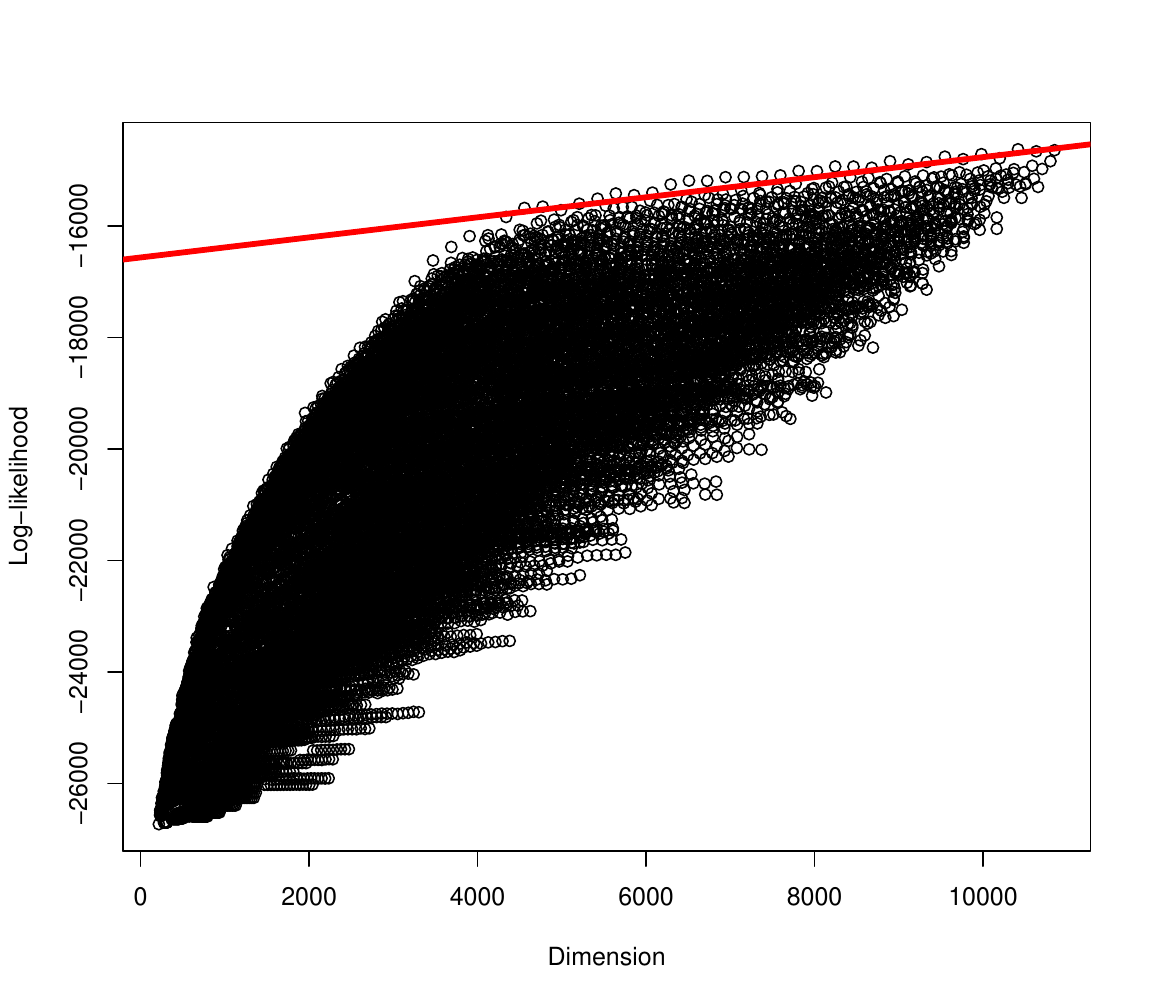}
\else
\includegraphics[width=2.00in,height=2.00in]{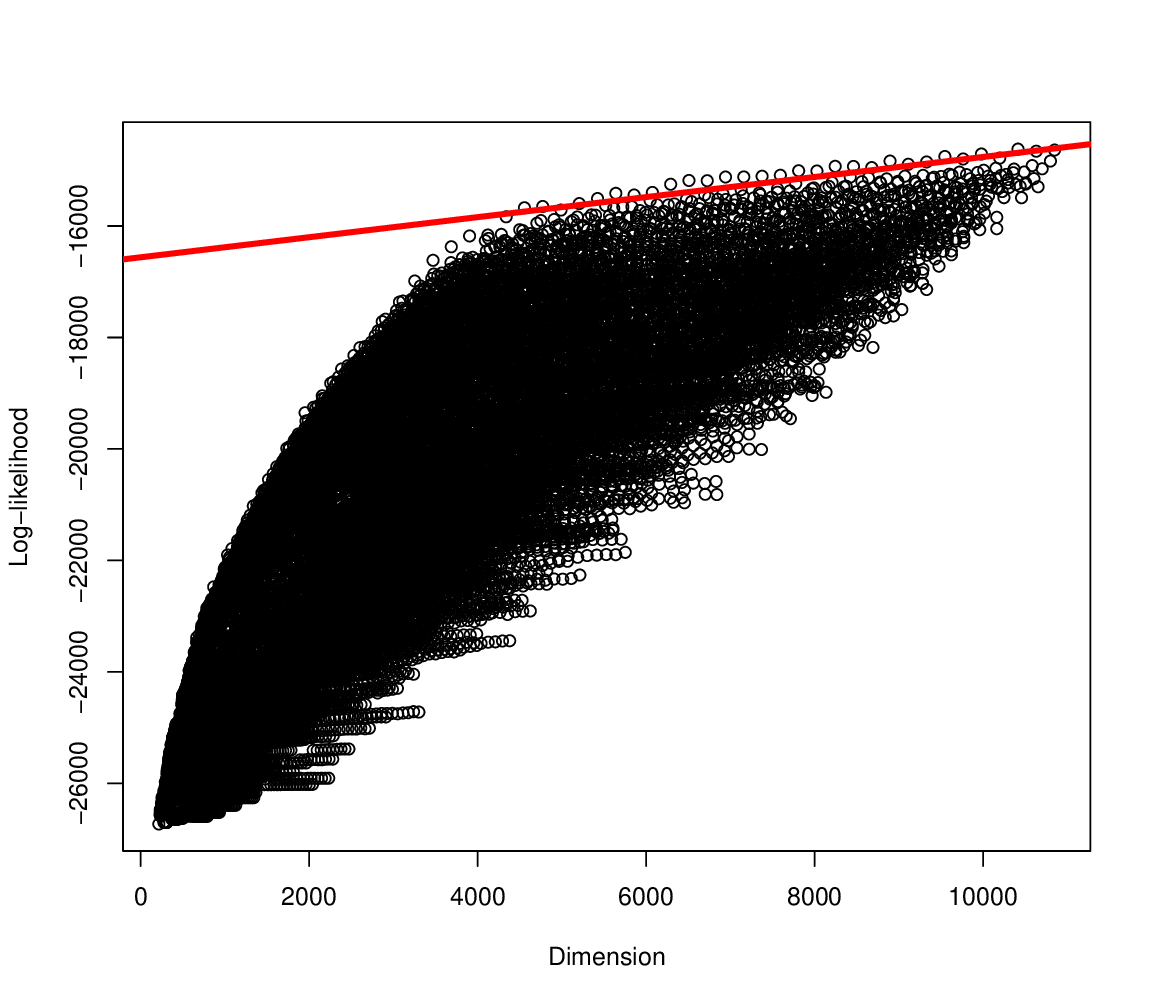}
\fi
}
\caption{Summary of the penalty calibration for chicken genotype data model
selection.}\label{Fig:ModelselectionChickenData}
\end{figure}

%s7 ###
\section{Conclusion}

We were able to simultaneously select variables and detect the number
of populations in the specific framework of multivariate multinomial
mixtures in our investigation of model selection via penalization. This
led to secondary clustering. Our main result provides an oracle
inequality, conditional on some lower bound on the penalty function.
The weakness of such a result is that the associated penalized
criterion is not directly usable. Nevertheless, it suggests a shape of
the penalty function, which is of the form $\pen_n(m)=\lambda D_m/n$,
where $\lambda=\lambda\left(n,\ \mathcal {C}\right)$ is a parameter
that is dependent on the data and on the collection of the competing
models. In practice, $\lambda$ is calibrated via slope
heuristics.\looseness=1

In our simulated experiments, the new criterion with penalty calibration
showed good behaviors regarding the density estimation and the
selection of the true model. It also performed well for both large
and reasonably small numbers of individuals. We are therefore able
to answer the question ``Which criterion for with sample size?''

The model dimension grew very rapidly in our modeling scenario. In
real experiments, the number of individuals may be small, and different
models with reduced dimensions may be necessary. %We currently work on
%models which cluster the populations differently for each variable, as
%well as models which allocate the same probability to several
%categories.
Possible models include those that cluster populations differently
for each variable, as well as models that allocate the same probability
to several categories in some clusters.

\appendix

%s8 ###
\section{Metric entropy with bracketing} \label{sec:entropy}

\begin{sloppypar}
We first provide a number of results concerning entropy with bracketing
that serve to prove Proposition~\ref{prop:entropy}. These results are
mainly adapted from \cite{GenoveseWasserman2000}, but several of them
were improved or rewritten in a more general form. These lemmas can be
regarded as a toolbox to calculate the metric entropy with bracketing
of complex models using the metric entropy of simpler elements.
\end{sloppypar}

We consider a measurable space $(A, {\mathcal A})$ and $\mu$ as a
$\sigma$-finite positive measure on $A$. We consider a model $\mathcal
{M}$, which is a set of probability density functions with respect to
$\mu$. All functions considered in the following are positive functions
in $\LL^1(\mu)$.

\begin{lem} \label{lem:maj-u-in-bracket}
Let $\varepsilon>0$. Let $[l,u]$ be a bracket in $\LL^1(\mu)$ with an
$\h$-diameter less than $\varepsilon$ and containing $s$, which is a
probability density function with respect to~$\mu$. Then
\[
\int l\, \ud\mu\leq1 \leq\int u\, \ud\mu\leq(1+\varepsilon)^2.
\]
\end{lem}

\begin{proof}
Two inequalities are immediate from $l\leq s \leq u$. The latter uses
the triangle inequality in $\LL^2(\mu)$ and the definition of $\h$:
\begin{align*}
\int u\, \ud\mu&= \int\left(\sqrt{l} + \left(\sqrt{u} - \sqrt
{l}\right) \right)^2 \ud\mu\\
&\leq\left(\sqrt{\int l\, \ud\mu} + \h(u,l)\right)^2 \\
&\leq(1+\varepsilon)^2.\hspace*{85pt}\qedhere
\end{align*}
\end{proof}

\begin{lem}[Bracketing entropy of product densities] \label{lem:entropy-prod}
Let $n\geq2$, and consider a collection $(A_i, {\mathcal A}_i, \mu
_i)_{1\leq i \leq n}$ of measured spaces. For any $1\leq i \leq n$, let
$\mathcal{M}_i$ be a collection of probability density functions on
$A_i$. % fulfilling~\ref{cond:M}.
Consider the product model
\[
\mathcal{M} = \left\{ s=\otimes_{i=1}^{n} s_i; \forall1\leq i \leq n,
s_i \in\mathcal{M}_i \right\}.
\]
$\mathcal{M}$ contains density functions on $A=\prod_{i=1}^{n} A_i$
with respect to $\mu=\otimes_{i=1}^{n} \mu_i$.

%$\mathcal{M}$ fulfills~\ref{cond:M} and, for
For any sequence of positive numbers $(\delta_i)_{1\leq i\leq n}$, if
$\varepsilon\geq\prod_{i=1}^{n} (1+\delta_i) -1$, then
% for any $\delta>0$, if $\varepsilon\geq(1+\delta)^n-1$ then
%
\[
H_{[\cdot]}\left(\varepsilon, \mathcal{M}, \h\right) \leq\sum
_{i=1}^{n} H_{[\cdot]}\left(\delta_i, \mathcal{M}_i, \h\right).
\]
\end{lem}

\begin{sloppypar}
\begin{proof}
%Let us consider some $s=\otimes_{i=1}^{n} s_i$ in $\mathcal{M}$. For
%$1\leq i \leq n$, let $\mathcal{M}_i'$, $A_i'$ and a sequence
%$(t_{i,k})_{k\geq1}$ be such as needed for $\mathcal{M}_i$ to satisfy
%\ref{cond:M}. Then, with the choice $t_k = \otimes_{i=1}^{n} t_{i,k}$
%and $A' = \prod_{i=1}^{n} A_i'$,~\ref{cond:M} is true for $
%\mathcal{M}$ too.
%
Let $\delta> 0$. For any $1\leq i \leq n$, let $[l_i, u_i]$ be a
bracket containing $s_i$, with an $\h$-diameter less than $\delta_i$.
Let $l = \otimes_{i=1}^{n} l_i$ and $u = \otimes_{i=1}^{n} u_i$. Then,
$s$ belongs to the bracket $[l, u]$, and we can compute its $\h
$-diameter as follows:
\begin{align*}
\h(l,u) &= \sqrt{\int_{A} \bigg( \sum_{j=1}^{n} \bigg(\prod
_{i=1}^{j-1} \sqrt{l_i} \prod_{i=j}^{n} \sqrt{u_i} - \prod_{i=1}^{j}
\sqrt{l_i} \prod_{i=j+1}^{n} \sqrt{u_i} \bigg) \bigg)^2\,\ud\mu}
\\
&\leq\sum_{j=1}^{n} \prod_{i=1}^{j-1} \sqrt{\int_{A_i} l_i\,\ud
\mu_i}
\prod_{i=j+1}^{n} \sqrt{\int_{A_i} u_i\,\ud\mu_i} \ \h(l_j, u_j)
\\
%&\leq\sum_{j=1}^{n} \delta\, (1+\delta)^{n-j} = (1+\delta)^n -1
&\leq\sum_{j=1}^{n} \delta_j \prod_{i=j+1}^{n} (1+\delta_i) %\\ &
= \prod_{j=1}^{n} (1+\delta_j) -1
\end{align*}
thanks to the triangle inequality and Lemma~\ref{lem:maj-u-in-bracket}
(empty products equal $1$).

Let $\varepsilon\geq\prod_{i=1}^{n} (1+\delta_i) -1$. For any $1\leq i
\leq n$, consider a minimal covering of $\mathcal{M}_i$ with brackets
of $\h$-diameter less than $\delta_i$. The previous process allows us
to build a covering of $\mathcal{M}$ with brackets of $\h$-diameter
less than $\varepsilon$. Thus, the minimal cardinality of such a
covering satisfies
\begin{equation*}
N_{[\cdot]}\left(\varepsilon, \mathcal{M}, \h\right) \leq\prod
_{i=1}^{n} N_{[\cdot]}\left(\delta_i, \mathcal{M}_i, \h\right).\qedhere
\end{equation*}
\end{proof}
\end{sloppypar}

\begin{lem}[Bracketing entropy of mixture densities] \label
{lem:entropy-mixture} Let $n\geq2$, and for any $1\leq i \leq n$, let
$\mathcal{M}_i$ be a set of probability density functions, all on the
same measured space
$(A, {\mathcal A}, \mu)$. % and fulfilling~\ref{cond:M}.
Consider the set of all mixture densities
\[
\mathcal{M} = \left\{ \sum_{i=1}^{n} \pi_i s_i: \mathbf{\pi}=(\pi
_i)_{1\leq i \leq n} \in\Ss_{n-1}; \forall1\leq i \leq n, s_i \in
\mathcal{M}_i \right\}.
\]
Then %$\mathcal{M}$ fulfills~\ref{cond:M}, and
for any $\delta>0$, $\eta>0$ and $\varepsilon\geq\delta+ \eta+
\delta\eta$,
\[
H_{[\cdot]}\left(\varepsilon, \mathcal{M}, \h\right) \leq
H_{[\cdot
]}\left(\delta, \Ss_{n-1}, \h\right) + \sum_{i=1}^{n} H_{[\cdot
]}\left
(\eta, \mathcal{M}_i, \h\right).
\]
\end{lem}

\begin{sloppypar}
\begin{proof}
%First, let us note that $\Ss_{n-1}$ is separable for its usual
%topology. Then, checking that $\mathcal{M}$ fulfills~\ref{cond:M} is
%easy, and we do not explicit it.
%
%We do not develop either the proof of the last relation,
We did not develop the proof because it is identical to \cite[proof of
Theorem~2]{GenoveseWasserman2000}. However, by using our Lemma~\ref
{lem:maj-u-in-bracket} instead of \cite
[Lemma~3]{GenoveseWasserman2000}, we obtain
\begin{align*}
\h^2(l,u) &\leq\eta^2 \, (1+\delta)^2 + \delta^2 + 2 \eta\,
\delta\,
(1+\delta) \\
&\leq\varepsilon^2.\hspace*{85pt}\qedhere
\end{align*}
\end{proof}
\end{sloppypar}

The following result merely restates Lemma~2 from
\cite{GenoveseWasserman2000}:

\begin{lem}[Bracketing entropy of the simplex] \label{lem:entropy-simplex}
Let $n\geq2$ be an integer. Let $\mu$ be the counting measure on $\{ 1,
\ldots, n\}$. We identify any probability on $\{1, \ldots, n\}$ with
its density $s\in\Ss_{n-1}$ with respect to $\mu$. Then, if $0 < \delta
\leq1$,
\[
H_{[\cdot]}\left(\delta, \Ss_{n-1}, \h\right) \leq(n-1) \ln
\left(\frac
{1}{\delta}\right) + \frac{\ln2 + \ln(n+1) + n \ln(2\pi e)}{2}.
\]
\end{lem}

In addition, the metric entropy of the collection of all Hardy-Weinberg
genotype distributions for a given variable is required to examine
Case~\ref{config:double-multinomial}.
\begin{lem}[Bracketing entropy of Hardy-Weinberg genotype
distributions] \label{lem:entropy-genotypes} Suppose that, for some
variable $l$, there exist $A_l\geq2$ different states. Let $\Omega_l$
be the collection of all genotype distributions
following the Hardy-Weinberg model~(\ref{HW_modele}). Then %$\Omega_l$
%fulfills~\ref{cond:M}, and
for any $\varepsilon>0$ and $\delta\geq\varepsilon\, (2+\varepsilon)$,
\[
H_{[\cdot]}\left(\delta, \Omega_l, \h\right) \leq H_{[\cdot
]}\left
(\varepsilon, \Ss_{A_l-1}, \h\right).
\]
\end{lem}

\begin{proof}
(\ref{HW_modele}) permits to associate a parameter $\alpha=(\alpha_1,
\ldots, \alpha_{A_l}) \in\Ss_{A_l-1}$ with any density in $\Omega_l$.
More generally, for any $\alpha\in[0,1]^{A_l}$, we define a function
\[
d_\alpha(x) = \left(2-\1_{x_{1}=x_{2}}\right) \alpha_{x_{1}} \alpha
_{x_{2}}
\]
on the set of all genotypes $x=\{x^1, x^2\}$ on $A_l$ states. Consider
some $\varepsilon>0$ and $d_\alpha\in\Omega_l$. Let $[l, u]$ be some
bracket containing $\alpha$, with an $\h$-diameter less
than~$\varepsilon$. Then $d_\alpha$ belongs to the bracket $[d_l,
d_u]$. The following calculates its diameter using
Lemma~\ref{lem:maj-u-in-bracket}:
\begin{align*}
\h^2(d_l, d_u) &= \sum_{a=1}^{A_l} \left(u_a - l_a\right)^2 + \sum
_{1\leq a<b\leq A_l} \left(\sqrt{2u_a u_b} - \sqrt{2l_a l_b} \right)^2
\\
&\leq\sum_{a=1}^{A_l} \sum_{b=1}^{A_l} \left(\sqrt{u_a u_b} -
\sqrt
{u_a l_b} + \sqrt{u_a l_b} - \sqrt{l_a l_b} \right)^2 \\
&\leq\left( \sqrt{\sum_{a=1}^{A_l} u_a \sum_{b=1}^{A_l} \left
(\sqrt
{u_b} - \sqrt{l_b}\right)^2} + \sqrt{\sum_{a=1}^{A_l} \left(\sqrt{u_a}
- \sqrt{l_a}\right)^2 \sum_{b=1}^{A_l} l_b} \right)^2 \\
&\leq\left( (1+\varepsilon)\, \varepsilon+ \varepsilon\right)^2
\end{align*}
Thus, $\h(d_l, d_u) \leq\varepsilon\, (2+\varepsilon)$.
%
%Let $(\alpha^{(k)})_{k\geq1}$ a sequence of elements of $\Ss_{A_l-1}
%\cap\Q^{A_l}$, which tends to $\alpha$ for the usual topology as $k$
%tends to infinity. Then, for any genotype $x=\{x^1, x^2\}$, $\ln d_{
%\alpha^{(k)}}(x)$ tends to $\ln d_{\alpha}(x)$. Therefore $\Omega_l$
%fulfills~\ref{cond:M}.
\end{proof}

\begin{proof}[Proof of Proposition~\ref{prop:entropy}]
We built the proof for Case~\ref{config:double-multinomial}. Case~\ref
{config:simple-multinomial} is very similar although it includes the
following simplification: we directly obtain $\Ss_{A_l-1}$ instead of
$\Omega_l$.

Using (\ref{Melange_selection}), we observe that a probability
$P_{(K,S)}\big(\cdot|\theta\big)$ is the product of two terms: the
first is a mixture density associated to the variables in $S$, and the
second is a product density on $\bigotimes_{l\notin S} \Omega_l$
associated to the other variables. Let $\mathcal{M}$ denote the
collection of all mixtures of $K$ densities in $\bigotimes_{l\in S}
\Omega_l$.

We first address the non-clustering variables. Given Lemma~\ref
{lem:entropy-prod} and Lemma~\ref{lem:entropy-genotypes}, %$
%\bigotimes_{l\notin S} \Omega_l$ fulfills~\ref{cond:M}. For
for any $\varepsilon\in(0,1)$,
\begin{align*}
H_{[\cdot]}\left((1+\varepsilon)^{2(L-|S|)}-1, \bigotimes_{l\notin S}
\Omega_l, \h\right) &\leq\sum_{l\notin S} H_{[\cdot]}\left
(\varepsilon
(2+\varepsilon), \Omega_l, \h\right) \\
&\leq\sum_{l\notin S} H_{[\cdot]}\left(\varepsilon, \Ss_{A_l-1},
\h
\right).
\end{align*}

Correspondingly,
\[
H_{[\cdot]}\left((1+\varepsilon)^{2|S|}-1, \bigotimes_{l\in S}
\Omega
_l, \h\right) \leq\sum_{l\in S} H_{[\cdot]}\left(\varepsilon,
\Ss
_{A_l-1}, \h\right).
\]\eject\noindent
Applying Lemma~\ref{lem:entropy-mixture}, we obtain %$\mathcal{M}$
%fulfills~\ref{cond:M} and
%
\begin{multline*}
H_{[\cdot]}\left((1+\varepsilon)^{2|S|+1}-1, \mathcal{M}, \h\right
) \\
\leq\1_{K\geq2} H_{[\cdot]}\left(\varepsilon, \Ss_{K-1}, \h \right) + K
\sum_{l\in S} H_{[\cdot]}\left(\varepsilon, \Ss_{A_l-1}, \h \right).
\end{multline*}

Applied to $\mathcal{M}$ and $\bigotimes_{l\notin S} \Omega_l$,
Lemma~\ref{lem:entropy-prod} gives %$\mathcal{M}_{(K, S)}$ fulfills
%\ref{cond:M}, and
for any $\varepsilon\in(0,1)$,
\begin{multline*}
H_{[\cdot]}\left(\eta_L(\varepsilon), \mathcal{M}_{(K, S)}, \h \right)
\\
\leq\1_{K\geq2} H_{[\cdot]}\left(\varepsilon, \Ss_{K-1}, \h \right) + K
\sum_{l\in S} H_{[\cdot]}\left(\varepsilon, \Ss_{A_l-1}, \h \right) +
\sum_{l\notin S} H_{[\cdot]}\left(\varepsilon, \Ss_{A_l-1}, \h \right).
\end{multline*}

At this stage, only Lemma~\ref{lem:entropy-simplex} remains to be
applied and the constants must be computed.
% \begin{align*}
% H_{[\cdot]}\left(\eta_L(\varepsilon), \mathcal{M}_{(K, S)}, \h\right)
%&\leq\left(\ln\left(\frac{1}{\varepsilon}\right) + \frac{\ln(2\pi
%e)}{2}\right) D_{(K, S)} \\
% &\quad+ \frac{\ln(4\pi e)}{2} \left(\1_{K\geq2} + K |S| + L-|S|
%\right) \\
% &\quad+\frac{1}{2}\Bigg( \1_{K\geq2} \ln(K+1) + \sum_{l=1}^{L}
%\ln(A_l+1) \\
% &\quad+ (K-1) \sum_{l\in S} \ln(A_l+1)\Bigg).
% \end{align*}
\end{proof}

%s9 ###
\section{Establishing the penalty} \label{sec:phi}

First, some properties of function $\eta_L$ are necessary. Recall that
$\eta_L(\varepsilon) = (1+\varepsilon)^{L+1} - 1$ in Case~\ref
{config:simple-multinomial}, and $\eta_L(\varepsilon) = (1+\varepsilon
)^{2L+1} - 1$ in Case~\ref{config:double-multinomial}.

\begin{lem}[Properties of the function $\eta_L$] \label{lem:prop-eta}
We consider the function $\eta_L$ defined in Proposition~\ref
{prop:entropy}, from $\R_+$ into $\R_+$. The function $\eta_L$ is
nonnegative, increasing, and convex. $\eta_L(0)=0$, and $\eta
_L'(0)=L+1$ in Case~\ref{config:simple-multinomial}, whereas $\eta
_L'(0)=2L+1$ in Case~\ref{config:double-multinomial}.
\end{lem}

\begin{proof}[Proof of Proposition~\ref{prop:phi-m}]
Let $0<\sigma\leq\eta_L(1)$, and let $\delta= \eta_L^{-1}(\sigma)$.
Then, for any $u \in\mathcal{M}_{m}$,
\begin{align*}
&\int_{0}^{\sigma} \sqrt{H_{[\cdot]}\left(x, \mathcal
{M}_{m}(u,\sigma
), \h\right)} \ud x \\
&\qquad\leq\sum_{j=1}^{\infty} \int_{\eta_L(2^{-j}\delta)}^{\eta
_L(2^{-j+1}\delta)} \sqrt{H_{[\cdot]}\left(x, \mathcal{M}_{m}, \h
\right
)} \ud x \\
&\qquad\leq\sum_{j=1}^{\infty} \left(\eta_L\left(2^{-j+1}\delta
\right
) - \eta_L\left(2^{-j}\delta\right)\right) \sqrt{C_{m} - D_{m}
\ln
\delta+ D_{m} j\ln2} \\
&\qquad\leq\eta_L(\delta) \sqrt{C_{m} - D_{m} \ln\delta} \\
&\qquad\quad+ \sqrt{D_{m}\ln2}\, \sum_{j=1}^{\infty} \sqrt{j}\,
\left
(\eta_L\left(2^{-j+1}\delta\right) - \eta_L\left(2^{-j}\delta
\right
)\right).
\end{align*}
The last term of this sum is addressed in the following:
\begin{align*}
\sum_{j=1}^{\infty} \sqrt{j}\, \left(\eta_L\left(2^{-j+1}\delta
\right)
- \eta_L\left(2^{-j}\delta\right)\right)
&\leq\sum_{j=1}^{\infty} j\, \left(\eta_L\left(2^{-j+1}\delta
\right)
- \eta_L\left(2^{-j}\delta\right)\right) \\
%&= \sum_{k=1}^{\infty} \sum_{j=k}^{\infty} \left(\eta_L\left(2^{-j+1}
%\delta\right) - \eta_L\left(2^{-j}\delta\right)\right) \\
&= \sum_{k=1}^{\infty} \eta_L\left(2^{-k+1}\delta\right) \\
&\leq\sum_{k=1}^{\infty} 2^{-k+1} \eta_L(\delta) = 2 \sigma.
\end{align*}
So
\[
\int_{0}^{\sigma} \sqrt{H_{[\cdot]}\left(x, \mathcal
{M}_{m}(u,\sigma),
\h\right)} \ud x \leq\phi_m(\sigma).
\]
Because $\eta_L$ is increasing, $\phi_m(x)/x$ is decreasing. To verify
that $\phi_m$ is nondecreasing, it is sufficient to prove that the
function $f(x)= x \sqrt{b-\ln\eta_L^{-1}(x)}$ is nondecreasing on $(0,
\eta_L(1)]$, where $b = \frac{C_{m}}{D_{m}}$. From (\ref{eq:def-RKS}),
we get $C_{m} > \frac{\ln(2\pi e)}{2} D_{m} > D_{m}$, so that $b>1$.
Calculus gives
\[
f'(x) = \sqrt{b-\ln\eta_L^{-1}(x)} - \frac{x}{2 \eta_L^{-1}(x)\,
\eta
_L'\left(\eta_L^{-1}(x)\right) \sqrt{b-\ln\eta_L^{-1}(x)}}.
\]

Let $y\in(0, 1]$. The function $\eta_L$ is convex on $(0, 1]$, which
entails $\frac{\eta_L(y)}{y\, \eta_L'(y)} \leq1$. Thus
\begin{equation*}
\sqrt{b-\ln y}\, f'\left(\eta_L(y)\right) \geq b - \ln y -1/2
> 0.\qedhere
\end{equation*}
\end{proof}

%$\xi= \dfrac{4 \sqrt{A_{\max}}\sqrt{L}}{2^{L+1}-1}$ in Case~
%\ref{config:simple-multinomial} and $\xi= \dfrac{4 \sqrt{A_{\max}}
%\sqrt{L}}{2 (1+3\sqrt{2})^{L}-1}$ in Case~
%\ref{config:double-multinomial}

\begin{proof}[Proof of Lemma~\ref{lem:sigma-less-than-eta}]
For any $\sigma>0$ such that $\sigma> \frac{\phi_m(\sigma)}{\sqrt {n}\,
\sigma}$, we have $\sigma>\sigma_m$, because $x\mapsto\frac {\phi
_m(x)}{x}$ is a nonincreasing function. Therefore, to obtain $\sigma
_m=\frac{\phi_m(\sigma_m)}{\sqrt{n}\, \sigma_m}<\eta_L(1)$, it suffices
that $\sqrt{n} > \frac{\phi_m(\eta_L(1))}{\eta_L^2(1)}$.

For all $1\leq l \leq L$, $A_l\geq2$. Because $\frac{1}{2} \ln(1+x)
\leq x-1$ for $x\geq2$, we obtain the following bounds:
%
%e13 ###
\begin{equation} \label{eq:R-K-S}
\frac{1+\ln(2\pi)}{2} D_{m} \leq C_{m} \leq\left( 2+\ln(2\pi) +
\frac{\ln2}{2}\right) D_{m}.
\end{equation}
Conversely, we have
\begin{equation*} %\label{eq:D-KLM}
D_{m} \leq K\, L\, A_{\max}.
\end{equation*}
Therefore,
\begin{align*}
\frac{\phi_m(\eta_L(1))}{\eta_L^2(1)} &\leq\frac{\left(2\sqrt
{\ln
(2)}+\sqrt{2+\ln(2\pi)+\ln(2)/2}\right)\sqrt{D_m}}{\eta_L(1)}\\
& < \frac{4 \sqrt{D_{m}}}{\eta_L(1)} %\\&
< \frac{4\sqrt{KLA_{\max}}}{\eta_L(1)}%\\&
< \frac{4\sqrt{LA_{\max}}}{\eta_L(1)}\sqrt{n}.
\end{align*}
Recall that only models with $K\leq n$ were considered. Thus, we have
$\sigma_m < \eta_L(1)$ as soon as $\xi=\frac{4\sqrt{LA_{\max }}}{\eta
_L(1)}\leq1$.
\end{proof}

\begin{proof}[Proof of Lemma~\ref{lem:sum-x-m}]
We define $\delta= 1/2$, from which $e^{-x_m} = \delta^{D_{m}}$.
Considering the collection $\colmodels$, we can distinguish two cases:
$K=1$ and $S=\emptyset$, or $K\geq2$ and $S\neq\emptyset$. Thus, using
(\ref{eq:Dimension}),
\begin{align*}
\sum_{m\in\colmodels} e^{-x_m} &= \delta^{\sum_{l=1}^{L} \left
(A_l-1\right)} \Bigg(1+ \sum_{S\neq\emptyset} \sum_{K\geq2}
\left
(\delta^{1+\sum_{l\in S}\left(A_l-1\right)}\right)^{K-1} \Bigg) \\
&= \delta^{\sum_{l=1}^{L} \left(A_l-1\right)} \Bigg(1+ \sum
_{S\neq
\emptyset} \frac{\delta^{1+\sum_{l\in S}\left(A_l-1\right
)}}{1-\delta
^{1+\sum_{l\in S}\left(A_l-1\right)}} \Bigg) \\
&\leq\delta^L \Bigg(1+\frac{\delta}{1-\delta} \sum_{S\neq
\emptyset}
\delta^{|S|}\Bigg) \\
&= \delta^L (1+\delta)^L. \hspace*{85pt}\qedhere%\delta= 1/2
\end{align*}
\end{proof}

\begin{proof}[Proof of Lemma~\ref{lem:pen}]
The function $\eta_L^{-1}$ is nondecreasing and concave, and it is
given by
\begin{equation*}
\eta_L^{-1}\left(x\right)=\left\{
\begin{array}{ll}
(x+1)^{\frac{1}{L+1}}-1 & \text{in Case~\ref
{config:simple-multinomial},}\\
(x+1)^{\frac{1}{2L+1}}-1 & \text{in Case~\ref{config:double-multinomial}.}
\end{array}
\right.
\end{equation*}
For any $0\leq x \leq\eta_L(1)$,
\[
\eta_L^{-1}(x) \geq\frac{\eta_L^{-1}(2)}{2}\, \min(x, 2).
\]

However, using (\ref{eq:sigma-tilde}) and (\ref{eq:R-K-S}), we obtain
%
%e14 ###
\begin{equation} \label{eq:bounds-sigma-tilde}
\widetilde{\sigma}_{m} \geq C_1 \sqrt{\frac{D_{m}}{n}} \geq C_1
\sqrt
{\frac{L}{n}},
\end{equation}
where $C_1 = 2\sqrt{\ln2} + \sqrt{\frac{1+\ln(2\pi)}{2}} >2\sqrt{2}$.
Therefore,
\[
-\ln\eta_L^{-1}\left(\widetilde{\sigma}_{m}\right) \leq-\ln
\left
(\frac{\eta_L^{-1}(2)}{2}\right) -\ln2 + \max\left(0, \frac
{1}{2}\left
(\ln n - \ln L - \ln2\right) \right).
\]

Consider Case~\ref{config:simple-multinomial}. Because $\eta_L$ is a
convex function and $\eta_L'(0)=L+1$,
\[
\eta_L^{-1}(2) \leq\frac{2}{L+1}.
\]
Then,
\[
\eta_L\left(\frac{2}{L+1}\right) = \left(1+ \frac{2}{L+1}\right)^{L+1}
-1 \leq e^2-1.
\]
Therefore,
\[
\frac{\eta_L^{-1}(2)}{2} \geq\frac{2/(L+1)}{\eta_L\left
(2/(L+1)\right
)} \geq\frac{2}{(e^2-1)(L+1)}.
\]

\begin{sloppypar}
Then, considering Case~\ref{config:double-multinomial} in the same
manner, $\eta_L^{-1}(2) \leq\frac{2}{2L+1}$, $\eta_L(\frac
{2}{2L+1})\leq e^2-1$. This leads to
\[
\frac{\eta_L^{-1}(2)}{2} \geq\frac{2}{(e^2-1)(2L+1)},
\]
which is valid in both cases.
\end{sloppypar}

Therefore,
\[
-\ln\left(\frac{\eta_L^{-1}(2)}{2}\right) \leq\ln(e^2-1) + \ln L +
\ln(5/4)
\]
and
\[
-\ln\eta_L^{-1}\left(\widetilde{\sigma}_{m}\right) \leq\ln(e^2-1)
-\frac{7}{2} \ln2 + \ln5 + \max\left(\frac{1}{2} \ln n + \frac{1}{2}
\ln L, \frac{\ln2}{2} + \ln L\right).
\]

Using (\ref{eq:sigma-upper}), we obtain
\begin{align*}
\sigma_m^2+\frac{x_m}{n} &\leq\frac{D_{m}}{n}\, \left(\left
(2\sqrt
{\ln2} + \sqrt{2+\ln(2\pi) + \frac{\ln2}{2}-\ln\eta_L^{-1}\left
(\widetilde{\sigma}_{m}\right)} \right)^2 + \ln2\right) \\
&\leq\frac{D_{m}}{n}\, \Bigg(3\sqrt{\ln2} + \sqrt{2+\ln(2\pi) - 3
\ln2 + \ln5 + \ln(e^2-1)} \\
&\qquad+ \sqrt{\max\left( \frac{\ln n + \ln L}{2}, \frac{\ln
2}{2} +
\ln L\right)} \Bigg)^2 \\
&\leq\frac{D_{m}}{n}\, \left(5 + \sqrt{\max\left( \frac{\ln n +
\ln
L}{2}, \frac{\ln2}{2} + \ln L\right)} \right)^2,
\end{align*}
which is the desired result.
\end{proof}

\section*{Acknowledgments}

The authors would like to thank Elisabeth Gassiat, Pascal Massart, and
Gilles Celeux for their comments and advice. Thank you to Nathalie
Akakpo, Nicolas Verzelen, and Cathy Maugis whose helpful discussions we
very much appreciated.

\bibliographystyle{imsart-nameyear}%apalike
%\bibliography{PenaltyCalibration}

\begin{thebibliography}{27}
% BibTex style file: imsart-nameyear.bst, 2013-01-28
% Default style options (sort=1,type=nameyear).
% Used options (sort=1,type=nameyear).

%
%
%b1 ###
\bibitem[\protect\citeauthoryear{Arlot and Massart}{2009}]{ArlotMassart2009}
\begin{barticle}[author]
\bauthor{\bsnm{Arlot},~\bfnm{Sylvain}\binits{S.}} \AND
\bauthor{\bsnm{Massart},~\bfnm{Pascal}\binits{P.}}
(\byear{2009}).
\btitle{Data-driven calibration of penalties for least-squares regression}.
\bjournal{J. Mach. Learn. Res.}
\bvolume{10}
\bpages{245--279}.
\end{barticle}
\endbibitem

%
%
%b2 ###
\bibitem[\protect\citeauthoryear{Asuncion and
Newman}{2007}]{Asuncion+Newman:2007}
\begin{bmisc}[author]
\bauthor{\bsnm{Asuncion},~\bfnm{A.}\binits{A.}} \AND
\bauthor{\bsnm{Newman},~\bfnm{D.~J.}\binits{D.~J.}}
(\byear{2007}).
\btitle{{UCI} Machine Learning Repository}.
\end{bmisc}
\endbibitem

%
%
%b3 ###
\bibitem[\protect\citeauthoryear{Bai, Rao and Wu}{1999}]{BaiRaoWu1999}
\begin{barticle}[author]
\bauthor{\bsnm{Bai},~\bfnm{Z.}\binits{Z.}},
\bauthor{\bsnm{Rao},~\bfnm{C.~R.}\binits{C.~R.}} \AND
\bauthor{\bsnm{Wu},~\bfnm{Y.}\binits{Y.}}
(\byear{1999}).
\btitle{Model selection with data-oriented penalty}.
\bjournal{J. Statist. Plann. Inference}
\bvolume{77}
\bpages{102--117}.
\end{barticle}
\MR{1677811}
\endbibitem

%
%
%b4 ###
\bibitem[\protect\citeauthoryear{Biernacki, Celeux and
Govaert}{2000}]{biernacki2000}
\begin{barticle}[author]
\bauthor{\bsnm{Biernacki},~\bfnm{C.}\binits{C.}},
\bauthor{\bsnm{Celeux},~\bfnm{G.}\binits{G.}} \AND
\bauthor{\bsnm{Govaert},~\bfnm{G.}\binits{G.}}
(\byear{2000}).
\btitle{Assessing a mixture model for clustering with the integrated completed
likelihood}.
\bjournal{{IEEE} Trans. Pattern Anal.}
\bvolume{22}
\bpages{719--725}.
\end{barticle}
\endbibitem

%
%
%b5 ###
\bibitem[\protect\citeauthoryear{Birg{\'e} and
Massart}{2007}]{BirgeMassart2007}
\begin{barticle}[author]
\bauthor{\bsnm{Birg{\'e}},~\bfnm{Lucien}\binits{L.}} \AND
\bauthor{\bsnm{Massart},~\bfnm{Pascal}\binits{P.}}
(\byear{2007}).
\btitle{Minimal penalties for {G}aussian model selection}.
\bjournal{Probab. Theory Related Fields}
\bvolume{138}
\bpages{33--73}.
\end{barticle}
\MR{2288064}
\endbibitem

%
%
%b6 ###
\bibitem[\protect\citeauthoryear{Celeux and Govaert}{1991}]{CeleuxGovaert91}
\begin{barticle}[author]
\bauthor{\bsnm{Celeux},~\bfnm{Gilles}\binits{G.}} \AND
\bauthor{\bsnm{Govaert},~\bfnm{G{\'e}rard}\binits{G.}}
(\byear{1991}).
\btitle{Clustering criteria for discrete data and latent class models}.
\bjournal{J. Classif.}
\bvolume{8}
\bpages{157--176}.
\bdoi{10.1007/BF02616237}
\end{barticle}
\endbibitem

%
%
%b7 ###
\bibitem[\protect\citeauthoryear{Celeux, Hurn and
Robert}{2000}]{Celeux99computational}
\begin{barticle}[author]
\bauthor{\bsnm{Celeux},~\bfnm{Gilles}\binits{G.}},
\bauthor{\bsnm{Hurn},~\bfnm{Merrilee}\binits{M.}} \AND
\bauthor{\bsnm{Robert},~\bfnm{Christian~P.}\binits{C.~P.}}
(\byear{2000}).
\btitle{Computational and inferential difficulties with mixture posterior
distributions}.
\bjournal{J. Am. Stat. Assoc.}
\bvolume{95}
\bpages{957--970}.
\end{barticle}
\MR{1804450}
\endbibitem

%
%
%b8 ###
\bibitem[\protect\citeauthoryear{Chen, Forbes and Francois}{2006}]{Chen2006}
\begin{barticle}[author]
\bauthor{\bsnm{Chen},~\bfnm{C.}\binits{C.}},
\bauthor{\bsnm{Forbes},~\bfnm{F.}\binits{F.}} \AND
\bauthor{\bsnm{Francois},~\bfnm{O.}\binits{O.}}
(\byear{2006}).
\btitle{Fastruct: Model-based clustering made faster}.
\bjournal{Molecular Ecology Notes}
\bvolume{6}
\bpages{980--983}.
\end{barticle}
\endbibitem

%
%
%b9 ###
\bibitem[\protect\citeauthoryear{Collins and Lanza}{2010}]{CollinsLanza}
\begin{bbook}[author]
\bauthor{\bsnm{Collins},~\bfnm{Linda~M.}\binits{L.~M.}} \AND
\bauthor{\bsnm{Lanza},~\bfnm{Stephanie~T.}\binits{S.~T.}}
(\byear{2010}).
\btitle{Latent Class and Latent Transition Analysis: With Applications
in the
Social, Behavioral, and Health Sciences}.
\bseries{Wiley Series in Probability and Statistics}.
\bpublisher{Wiley}.
\end{bbook}
\endbibitem

%
%
%b10 ###
\bibitem[\protect\citeauthoryear{Corander et~al.}{2008}]{Corander2008}
\begin{barticle}[author]
\bauthor{\bsnm{Corander},~\bfnm{Jukka}\binits{J.}},
\bauthor{\bsnm{Marttinen},~\bfnm{Pekka}\binits{P.}},
\bauthor{\bsnm{Sir{\'e}n},~\bfnm{Jukka}\binits{J.}} \AND
\bauthor{\bsnm{Tang},~\bfnm{Jing}\binits{J.}}
(\byear{2008}).
\btitle{Enhanced {B}ayesian modelling in {BAPS} software for learning genetic
structures of populations}.
\bjournal{BMC Bioinformatics}
\bvolume{9}
\bpages{539}.
\end{barticle}
\endbibitem

%
%
%b11 ###
\bibitem[\protect\citeauthoryear{Dempster, Lairdsand and
Rubin}{1977}]{Dempster1977}
\begin{barticle}[author]
\bauthor{\bsnm{Dempster},~\bfnm{A.~P.}\binits{A.~P.}},
\bauthor{\bsnm{Lairdsand},~\bfnm{N.~M.}\binits{N.~M.}} \AND
\bauthor{\bsnm{Rubin},~\bfnm{D.~B.}\binits{D.~B.}}
(\byear{1977}).
\btitle{Maximum likelihood from incomplete data via the {EM} algorithm}.
\bjournal{J. Royal Statist. Soc. Series B}
\bvolume{39}
\bpages{1--38}.
\end{barticle}
\MR{0501537}
\endbibitem

%
%
%b12 ###
\bibitem[\protect\citeauthoryear{Genoveve and
Wasserman}{2000}]{GenoveseWasserman2000}
\begin{barticle}[author]
\bauthor{\bsnm{Genoveve},~\bfnm{C.~R.}\binits{C.~R.}} \AND
\bauthor{\bsnm{Wasserman},~\bfnm{Larry}\binits{L.}}
(\byear{2000}).
\btitle{Rates of convergence for the {G}aussian mixture sieve}.
\bjournal{Ann. Statist.}
\bvolume{28}
\bpages{1105--1127}.
\end{barticle}
\MR{1810921}
\endbibitem

%
%
%b13 ###
\bibitem[\protect\citeauthoryear{Goodman}{1974}]{Goodman74}
\begin{barticle}[author]
\bauthor{\bsnm{Goodman},~\bfnm{Leo~A.}\binits{L.~A.}}
(\byear{1974}).
\btitle{Exploratory latent structure analysis using both identifiable and
unidentifiable models}.
\bjournal{Biometrika}
\bvolume{61}
\bpages{215--231}.
\end{barticle}
\MR{0370936}
\endbibitem

%
%
%b14 ###
\bibitem[\protect\citeauthoryear{Latch et~al.}{2006}]{Latch2006}
\begin{barticle}[author]
\bauthor{\bsnm{Latch},~\bfnm{E.~K.}\binits{E.~K.}},
\bauthor{\bsnm{Dharmarajan},~\bfnm{Guha}\binits{G.}},
\bauthor{\bsnm{Glaubitz},~\bfnm{Jeffrey~C.}\binits{J.~C.}} \AND
\bauthor{\bsnm{Rhodes},~\bfnm{Olin~E.}\binits{O.~E.}~\bsuffix{Jr.}}
(\byear{2006}). \btitle{Relative performance of {B}ayesian clustering
software for inferring population substructure and individual
assignment at low levels of population differentiation}.
\bjournal{Conservation Genetics} \bvolume{7} \bpages{295}.
\end{barticle}
\endbibitem

%
%
%b15 ###
\bibitem[\protect\citeauthoryear{Lebarbier}{2002}]{LebarbierPhDThesis}
\begin{bphdthesis}[author]
\bauthor{\bsnm{Lebarbier},~\bfnm{{\'E}milie}\binits{{\'E}.}}
(\byear{2002}).
\btitle{Quelques approches pour la d{\'e}tection de rupture {\`a}
horizon fini}
\btype{PhD thesis},
\bpublisher{Univ Paris-Sud},
\baddress{F-91405 Orsay}.
\end{bphdthesis}
\endbibitem

%
%
%b16 ###
\bibitem[\protect\citeauthoryear{Massart}{2007}]{Massart2007}
\begin{bbook}[author]
\bauthor{\bsnm{Massart},~\bfnm{Pascal}\binits{P.}}
(\byear{2007}).
\btitle{Concentration inequalities and model selection}.
\bseries{Lecture Notes in Mathematics}
\bvolume{1896}.
\bpublisher{Springer-Verlag}, \baddress{Berlin}.
\end{bbook}
\MR{2319879}
\endbibitem

%
%
%b17 ###
\bibitem[\protect\citeauthoryear{Maugis and
Michel}{2011a}]{MaugisMichel2008thm}
\begin{barticle}[author]
\bauthor{\bsnm{Maugis},~\bfnm{Cathy}\binits{C.}} \AND
\bauthor{\bsnm{Michel},~\bfnm{Bertrand}\binits{B.}}
(\byear{2011}a).
\btitle{A non asymptotic penalized criterion for {G}aussian mixture model
selection}.
\bjournal{ESAIM: P\&S}
\bvolume{15}
\bpages{41--68}.
\end{barticle}
\MR{2870505}
\endbibitem

%
%
%b18 ###
\bibitem[\protect\citeauthoryear{Maugis and
Michel}{2011b}]{MaugisMichel2008slope}
\begin{barticle}[author]
\bauthor{\bsnm{Maugis},~\bfnm{Cathy}\binits{C.}} \AND
\bauthor{\bsnm{Michel},~\bfnm{Bertrand}\binits{B.}} (\byear{2011}b).
\btitle{Data-driven penalty calibration: A~case study for Gaussian
mixture model selection}. \bjournal{ESAIM: P\&S} \bvolume{15}
\bpages{320--339}.
\end{barticle}
\MR{2870518}
\endbibitem

%
%
%b19 ###
\bibitem[\protect\citeauthoryear{McCutcheon}{1987}]{McCutcheon}
\begin{bbook}[author]
\bauthor{\bsnm{McCutcheon},~\bfnm{A.~L.}\binits{A.~L.}}
(\byear{1987}).
\btitle{Latent Class Analysis}.
\bseries{Quantitative Applications in the Social Sciences}
\bvolume{64}.
\bpublisher{Sage Publications}, \baddress{Thousand Oaks, California}.
\end{bbook}
\endbibitem

%
%
%b20 ###
\bibitem[\protect\citeauthoryear{McLachlan and Peel}{2000}]{McLachlanPeel}
\begin{bbook}[author]
\bauthor{\bsnm{McLachlan},~\bfnm{Geoffrey}\binits{G.}} \AND
\bauthor{\bsnm{Peel},~\bfnm{David}\binits{D.}}
(\byear{2000}).
\btitle{Finite Mixture Models}.
\bseries{Wiley Series in Probability and Statistics}.
\bpublisher{Wiley}.
\end{bbook}
\MR{1789474}
\endbibitem

%
%
%b21 ###
\bibitem[\protect\citeauthoryear{Nadif and Govaert}{1998}]{Nadif1998}
\begin{barticle}[author]
\bauthor{\bsnm{Nadif},~\bfnm{M.}\binits{M.}} \AND
\bauthor{\bsnm{Govaert},~\bfnm{G{\'e}rard}\binits{G.}} (\byear{1998}).
\btitle{Clustering for binary data and mixture models -- choice of the
model}. \bjournal{Appl. Stoch. Models Data Anal.} \bvolume{13}
\bpages{269--278}.
\end{barticle}
\endbibitem

%
%
%b22 ###
\bibitem[\protect\citeauthoryear{Pritchard, Stephens and
Donnelly}{2000}]{Pritchard2000}
\begin{barticle}[author]
\bauthor{\bsnm{Pritchard},~\bfnm{J.~K.}\binits{J.~K.}},
\bauthor{\bsnm{Stephens},~\bfnm{M.}\binits{M.}} \AND
\bauthor{\bsnm{Donnelly},~\bfnm{P.}\binits{P.}}
(\byear{2000}).
\btitle{{I}nference of population structure using multilocus genotype data.}
\bjournal{Genetics}
\bvolume{155}
\bpages{945--59}.
\end{barticle}
\endbibitem

%
%
%b23 ###
\bibitem[\protect\citeauthoryear{Rigouste, Capp{\'e} and
Yvon}{2006}]{Rigouste_Cappe_Yvon_2006}
\begin{barticle}[author]
\bauthor{\bsnm{Rigouste},~\bfnm{Lo{\"i}s}\binits{L.}},
\bauthor{\bsnm{Capp{\'e}},~\bfnm{Olivier}\binits{O.}} \AND
\bauthor{\bsnm{Yvon},~\bfnm{Fran\c{c}ois}\binits{F.}}
(\byear{2006}).
\btitle{Inference and evaluation of the multinomial mixture model for text
clustering}.
\bjournal{Inform. Process. Manag.}
\bvolume{43}
\bpages{1260--1280}.
\end{barticle}
\endbibitem

%
%
%b24 ###
\bibitem[\protect\citeauthoryear{Rosenberg et~al.}{2001}]{Rosenberg2001a}
\begin{barticle}[author]
\bauthor{\bsnm{Rosenberg},~\bfnm{Noah~A.}\binits{N.~A.}},
\bauthor{\bsnm{Burke},~\bfnm{Terry}\binits{T.}},
\bauthor{\bsnm{Elo},~\bfnm{Kari}\binits{K.}},
\bauthor{\bsnm{Feldman},~\bfnm{Marcus~W.}\binits{M.~W.}},
\bauthor{\bsnm{Freidlin},~\bfnm{Paul~J.}\binits{P.~J.}},
\bauthor{\bsnm{Groenen},~\bfnm{Martien A.~M.}\binits{M.~A.~M.}},
\bauthor{\bsnm{Hillel},~\bfnm{Jossi}\binits{J.}},
\bauthor{\bsnm{Ma},~\bfnm{Asko}\binits{A.}},
\bauthor{\bsnm{Vignal},~\bfnm{Alain}\binits{A.}},
\bauthor{\bsnm{Wimmers},~\bfnm{Klaus}\binits{K.}} \AND
\bauthor{\bsnm{Weigend},~\bfnm{Steffen}\binits{S.}}
(\byear{2001}).
\btitle{Empirical evaluation of genetic clustering methods using multilocus
genotypes from 20 chicken breeds}.
\bjournal{Biotechnology}.
\end{barticle}
\endbibitem

%
%
%b25 ###
\bibitem[\protect\citeauthoryear{Toussile and
Gassiat}{2009}]{ToussileGassiat2009}
\begin{barticle}[author]
\bauthor{\bsnm{Toussile},~\bfnm{Wilson}\binits{W.}} \AND
\bauthor{\bsnm{Gassiat},~\bfnm{Elisabeth}\binits{E.}}
(\byear{2009}).
\btitle{Variable selection in model-based clustering using multilocus genotype
data}.
\bjournal{Adv. Data Anal. Classif.}
\bvolume{3}
\bpages{109--134}.
\end{barticle}
\MR{2551051}
\endbibitem

%
%
%b26 ###
\bibitem[\protect\citeauthoryear{Verzelen}{2009}]{Verzelen2009}
\begin{bphdthesis}[author]
\bauthor{\bsnm{Verzelen},~\bfnm{Nicolas}\binits{N.}}
(\byear{2009}).
\btitle{Adaptative estimation to regular {G}aussian {M}arkov random fields}
\btype{PhD thesis},
\bpublisher{Univ Paris-Sud}.
\end{bphdthesis}
\endbibitem

%
%
%b27 ###
\bibitem[\protect\citeauthoryear{Villers}{2007}]{Villers2007}
\begin{bphdthesis}[author]
\bauthor{\bsnm{Villers},~\bfnm{F.}\binits{F.}}
(\byear{2007}).
\btitle{Tests et selection de mod{\`e}les pour l'analyse de donn{\'e}es
prot{\'e}omiques et transcriptomiques}
\btype{PhD thesis},
\bpublisher{Univ Paris-Sud}.
\end{bphdthesis}
\endbibitem

\end{thebibliography}

\end{document}